# LARGE SAMPLE THEORY OF INTRINSIC AND EXTRINSIC SAMPLE MEANS ON MANIFOLDS—II

By Rabi Bhattacharya[1] and Vic Patrangenaru[2]

*University of Arizona and Texas Tech University*

This article develops nonparametric inference procedures for estimation and testing problems for means on manifolds. A central limit theorem for Fréchet sample means is derived leading to an asymptotic distribution theory of intrinsic sample means on Riemannian manifolds. Central limit theorems are also obtained for extrinsic sample means w.r.t. an arbitrary embedding of a differentiable manifold in a Euclidean space. Bootstrap methods particularly suitable for these problems are presented. Applications are given to distributions on the sphere $S^d$ (directional spaces), real projective space $\mathbb{R}P^{N-1}$ (axial spaces), complex projective space $\mathbb{C}P^{k-2}$ (planar shape spaces) w.r.t. Veronese–Whitney embeddings and a three-dimensional shape space $\Sigma_3^4$.

**1. Introduction.** Statistical inference for distributions on manifolds is now a broad discipline with wide-ranging applications. Its study has gained momentum in recent years, especially due to applications in biosciences and medicine, and in image analysis. Including in the substantial body of literature in this field are the books by Bookstein [10], Dryden and Mardia [15], Kendall, Barden, Carne and Le [33], Mardia and Jupp [41], Small [49] and Watson [52]. While much of this literature focuses on parametric or semi-parametric models, the present article aims at providing a general framework for nonparametric inference for location. This is a continuation of our earlier work [7, 8] where some general properties of extrinsic and intrinsic mean sets on general manifolds were derived, and the problem of consistency of the corresponding sample indices was explored. The main focus of the present article is the derivation of asymptotic distributions of intrinsic and extrinsic

Received June 2002; revised March 2004.
[1]Supported by NSF Grant DMS-04-06143.
[2]Supported by NSA Grant MDA 904-02-1-0082 and NSF Grant OMS-04-06151.

sample means and confidence regions based on them. We provide classical CLT-based confidence regions and tests based on them, as well as those based on Efron's bootstrap [17].

Measures of location and dispersion for distributions on a manifold $M$ were studied in [7, 8] as Fréchet parameters associated with two types of distances on $M$. If $j:M \to \mathbb{R}^k$ is an embedding, the Euclidean distance restricted to $j(M)$ yields the *extrinsic mean set* and the *extrinsic total variance*. On the other hand, a *Riemannian distance* on $M$ yields the *intrinsic mean set* and *intrinsic total variance*.

Recall that the *Fréchet mean* of a probability measure $Q$ on a complete metric space $(M, \rho)$ is the minimizer of the function $F(x) = \int \rho^2(x,y) Q(dy)$, when such a minimizer exists and is unique [21]. In general the set of minimizers of $F$ is called the *Fréchet mean set*. The intrinsic mean $\mu_I(Q)$ is the Fréchet mean of a probability measure $Q$ on a *complete d-dimensional* Riemannian manifold $M$ endowed with the geodesic distance $d_g$ determined by the Riemannian structure $g$ on $M$. It is known that if $Q$ is sufficiently concentrated, then $\mu_I(Q)$ exists [see Theorem 2.2(a)]. The *extrinsic mean* $\mu_E(Q) = \mu_{j,E}(Q)$ of a probability measure $Q$ on a manifold $M$ w.r.t. an embedding $j: M \to \mathbb{R}^k$ is the Fréchet mean associated with the restriction to $j(M)$ of the Euclidean distance in $\mathbb{R}^k$. In [8] it was shown that the extrinsic mean of $Q$ exists if the ordinary mean of $j(Q)$ is a *nonfocal point* of $j(M)$, that is, if there is a *unique* point $x_0$ on $j(M)$ having the smallest distance from the mean of $j(Q)$. In this case $\mu_{j,E}(Q) = j^{-1}(x_0)$.

It is easier to compute the intrinsic mean if the Riemannian manifold has zero curvature in a neighborhood containing supp $Q$ [45]. In particular this is the case for distributions on linear projective shape spaces [42]. If the manifold has nonzero curvature around supp $Q$, it is easier to compute the extrinsic sample mean. It may be pointed out that if $Q$ is highly concentrated as in our medical examples in [8] and in Section 5, the intrinsic and extrinsic means are virtually indistinguishable.

We now provide a summary of the main results in this article. Section 2 is devoted to nonparametric inference for the Fréchet mean of a probability measure $Q$ on a manifold $M$ for which there is a domain $U$ of a chart $\phi: U \to \mathbb{R}^d$ such that $Q(U) = 1$. In Theorem 2.1 it is shown that in this case, under some rather general assumptions, the image of the Fréchet sample mean under $\phi$ is asymptotically normally distributed around the image of the Fréchet mean of $Q$. This leads to the asymptotic distribution theory of the intrinsic sample mean on a Riemannian manifold $M$ (Theorems 2.2, 2.3). In Corollaries 2.3 and 2.4 bootstrap confidence regions are derived for the Fréchet mean, with or without a pivot.

Section 3 is devoted to asymptotics of extrinsic sample means. The ideas behind the main result here are essentially due to Hendriks and Landsman [27] and Patrangenaru [44]. The two approaches are somewhat different.

We present in this article an extension of the latter approach. Extrinsic means are commonly used in directional, axial and shape statistics. In the particular case of directional data analysis, that is, when $M = S^{d-1}$ is the unit sphere in $\mathbb{R}^d$, Fisher, Hall, Jing and Wood [19] provided an approach for inference using computationally efficient bootstrapping which gets around the problem of increased dimensionality associated with the embedding of the manifold $M$ in a higher-dimensional Euclidean space. In Corollary 3.2 confidence regions are derived for the extrinsic mean $\mu_{j,E}(Q)$. Nonparametric bootstrap methods on abstract manifolds are also derived in this section (Theorem 3.2, Proposition 3.2).

If one assumes that $Q$ has a nonzero absolutely continuous component with respect to the volume measure on $M$, then from some results of Bhattacharya and Ghosh [6], Babu and Singh [1], Beran [2] and Hall [24, 25], one derives bootstrap-based confidence regions for $\mu_E(Q)$ with coverage error $O_p(n^{-2})$ (Theorem 3.4) (also see [5, 9]). One may also use the nonpivotal bootstrap to construct confidence regions based on the percentile method of Hall [25] for general $Q$ with a coverage error no more than $O_p(n^{-d/(d+1)})$, where $d$ is the dimension of the manifold (see Remark 2.4 and Proposition 3.2). This is particularly useful in those cases where the asymptotic dispersion matrix is difficult to compute.

Section 4 applies the preceding theory to (i) real projective spaces $\mathbb{R}P^{N-1}$—the *axial spaces*, and (ii) complex projective spaces $\mathbb{C}P^{k-2}$—the *shape spaces*. Another application to products of real projective spaces $(\mathbb{R}P^m)^{k-m-1}$, or the so-called *projective shape spaces*, will appear in [42].

As an application of Corollary 3.3, large sample confidence regions for mean axes are described in Corollary 4.2. A similar application to projective shape spaces, combining bootstrap methods for directional data from [3], appears in [42]. Other applications to axial spaces are given in Theorem 4.3 and Corollary 4.4, and to planar shape spaces in Theorem 4.5.

Finally in Section 5 we apply the results of Sections 2 and 4 to construct (1) a 95% large sample confidence region for the intrinsic mean location of the magnetic South Pole from a directional data set given in [20], (2) simultaneous confidence intervals for the affine coordinates of the extrinsic sample mean shape in a medical application and (3) a test for the difference between three-dimensional mean shapes in a glaucoma detection problem.

**2. A central limit theorem for Fréchet sample means and bootstrapping.** A $d$-dimensional *differentiable manifold* is a locally compact separable Hausdorff space $M$, together with an *atlas* $\mathcal{A}_M$ comprising a family of *charts* $(U_\alpha, \phi_\alpha)$ of open sets $U_\alpha$ covering $M$, and for each $\alpha$ a homeomorphism $\phi_\alpha$ of $U_\alpha$ onto an open subset of $\mathbb{R}^d$ for which the transition maps $\phi_\alpha \cdot \phi_\beta^{-1} : \phi_\beta(U_\alpha \cap U_\beta) \to \phi_\alpha(U_\alpha \cap U_\beta)$ are of class $\mathcal{C}^\infty$. The sets $U_\alpha$ are often called *coordinate neighborhoods*. One may show that a differentiable manifold is metrizable. We briefly recall some basic notion associated with Riemannian manifolds. For details the reader may refer to any standard text on differential geometry (e.g., [13, 26], or [38]). A *Riemannian metric* $g$ on a differentiable manifold $M$ is a $\mathcal{C}^\infty$ symmetric positive definite tensor field of type $(2,0)$, that is, a family of inner products $g_p = \langle \cdot, \cdot \rangle_p$ on the tangent spaces $T_p M, p \in M$, which is differentiable w.r.t. $p$. A *Riemannian manifold* $M$ is a connected differentiable manifold endowed with a Riemannian metric $g$. The distance $\rho_g$ induced by $g$ is called the *geodesic distance*. For $p, q \in M, \rho_g(p,q)$ is the infimum of lengths $\int_a^b \langle \dot{x}(t), \dot{x}(t) \rangle_{x(t)}^{1/2} dt$ of all $\mathcal{C}^1$-curves $x(t), a \leq t \leq b$, with $x(a) = p, x(b) = q$. The minimizer satisfies a variational equation whose solution is a *geodesic curve*. There is a unique such geodesic curve $t \to \gamma(t)$ for any initial point $\gamma(0) = p$ and initial tangent vector $\dot{\gamma}(0) = v$. A classical result of Hopf and Rinow states that $(M, \rho_g)$ is complete as a metric space if and only if $(M, g)$ is *geodesically complete* [i.e., every geodesic curve $\gamma(t)$ is defined for all $t \in \mathbb{R}$]. These two equivalent properties of completeness are in turn equivalent to a third property: *all closed bounded subsets of $(M, \rho_g)$ are compact* ([13], pages 146 and 147).

Given $q \in M$, the *exponential map* $\text{Exp}_q : U \to M$ is defined on an open neighborhood $U$ of $0 \in T_q M$ by the correspondence $v \to \gamma_v(1)$, where $\gamma_v(t)$ is the unique geodesic satisfying $\gamma(0) = q, \dot{\gamma}(0) = v$, provided $\gamma(t)$ extends at least to $t = 1$. Thus if $(M, g)$ is geodesically complete or, equivalently, $(M, \rho_g)$ is complete as a metric space, then $\text{Exp}_q$ is defined on all of $T_q M$. In this article, unless otherwise specified, all *Riemannian manifolds are assumed to be complete*.

Note that if $\gamma(0) = p$ and $\gamma(t)$ is a geodesic, it is generally not true that the geodesic distance between $p$ and $q = \gamma(t_1)$, say, is minimized by $\gamma(t), 0 \leq t \leq t_1$ (consider, e.g., the great circles on the sphere $S^2$ as geodesics). Let $t_0 = t_0(p)$ be the supremum of all $t_1 > 0$ for which this minimization holds. If $t_0 < \infty$, then $\gamma(t_0)$ is the *cut point of $p$ along $\gamma$*. The *cut locus* $C(p)$ of $p$ is the union of all cut points of $p$ along all geodesics $\gamma$ starting at $p$ [e.g., $C(p) = \{-p\}$ on $S^2$].

In this article we deal with both intrinsic and extrinsic means. Hence we will often consider a general distance $\rho$ on a differentiable manifold $M$, but assume that $(M, \rho)$ is complete as a metric space. We consider only those probability measures $Q$ on $M$ for which the Fréchet mean $\mu_\mathcal{F} = \mu_\mathcal{F}(Q)$ exists. Moreover *we assume that there is a chart $(U, \phi)$ such that $Q(U) = 1$, and $\mu_\mathcal{F} \in U$*.

REMARK 2.1. The assumption above on the existence of a chart $(U, \phi)$ such that $Q(U) = 1$ is less restrictive than it may seem. If $g$ is a Riemannian structure on $M$ and $Q$ is absolutely continuous w.r.t. the volume measure, then, for any given $p$, the complement $U$ of the cut locus $C(p)$ is the domain of definition of such a local coordinate system (the coordinate map being the inverse of $\text{Exp}_p$, the exponential map at $p$) (see [38], page 100, for details).

EXAMPLE 2.1. For the $d$-dimensional unit sphere, $M = S^d = \{p \in \mathbb{R}^{d+1} : \|p\| = 1\}$, with the Riemannian metric induced by the Euclidean metric on $\mathbb{R}^{d+1}$, the exponential map at a given point $p \in S^d$ is defined on the tangent space $T_p M$ and is given by

$$(2.1) \quad \text{Exp}_p(v) = \cos(\|v\|)p + \sin(\|v\|)\|v\|^{-1}v \qquad (v \in T_p S^d, v \neq 0).$$

If $x \in S^d, x \neq -p$, then there is a unique vector $u \in T_p M$ such that $x = \text{Exp}_p u$, and we will label this vector by $u = \text{Log}_p x$. Since $T_p S^d = \{v \in \mathbb{R}^{d+1}, v \cdot p = 0\}$, it follows that

$$(2.2) \quad \text{Log}_p x = (1 - (p \cdot x)^2)^{-1/2} \arccos(p \cdot x)(x - (p \cdot x)p).$$

In particular, for $d = 2$ we consider the orthobasis $e_1(p), e_2(p) \in T_p S^2$, where $p = (p_1, p_2, p_3)^t \in S^2 \setminus \{N, S\}$ $[N = (0, 0, 1), S = (0, 0, -1)]$:

$$(2.3) \quad \begin{aligned} e_1(p) &= ((p_1)^2 + (p_2)^2)^{-1/2}(-p_2, p_1, 0)^t, \\ e_2(p) &= (-((p_1)^2 + (p_2)^2)^{-1/2} p_1 p_3, \\ &\quad -(x^2 + y^2)^{-1/2} p_2 p_3, ((p_1)^2 + (p_2)^2)^{1/2})^t. \end{aligned}$$

The *logarithmic coordinates* of the point $x = (x_1, x_2, x_3)^T$ are given in this case by

$$(2.4) \quad \begin{aligned} u^1(p) &= e_1(p) \cdot \text{Log}_p x, \\ u^2(p) &= e_2(p) \cdot \text{Log}_p x. \end{aligned}$$

For computations one may use $a \cdot b = a^t b$.

Now the image measure $Q^\phi$ of $Q$ under $\phi$ has the Fréchet mean $\mu = \phi(\mu_\mathcal{F})$ w.r.t. the distance $\rho^\phi(u, v) := \rho(\phi^{-1}(u), \phi^{-1}(v)), u, v \in \phi(U)$. Similarly, if $X_i$ ($i = 1, \ldots, n$) are i.i.d. with common distribution $Q$ and defined on a probability space $(\Omega, \mathcal{A}, P)$, let $\mu_{n,\mathcal{F}}$ be a measurable selection from the Fréchet mean set (w.r.t. $\rho$) of the empirical $\hat{Q}_n = \frac{1}{n} \sum_{i=1}^n \delta_{X_i}$. Then $\mu_n = \phi(\mu_{n,\mathcal{F}})$ is a measurable selection from the Fréchet mean set (w.r.t. $\rho^\phi$) of $\hat{Q}_n^\phi = \frac{1}{n} \sum_{i=1}^n \delta_{\tilde{X}_i}$, where $\tilde{X}_i = \phi(X_i)$. Assuming twice continuous differentiability of $\theta \to (\rho^\phi)^2(u, \theta)$, write the Euclidean gradient as

$$(2.5) \quad \Psi(u; \theta) = \text{grad}_\theta (\rho^\phi)^2(u, \theta) = \left(\frac{\partial}{\partial \theta^r}(\rho^\phi)^2(u, \theta)\right)_{r=1}^d = (\Psi^r(u; \theta))_{r=1}^d.$$

Now $\mu$ is the point of minimum of

$$(2.6) \quad F^\phi(\theta) := \int (\rho^\phi)^2(u, \theta) Q^\phi(du)$$

and $\mu_n$ is a local minimum of

$$F_n^\phi(\theta) := \int (\rho^\phi)^2(u,\theta) \hat{Q}_n^\phi(du).$$

Therefore, one has the Taylor expansion

$$
\begin{aligned}
0 &= \frac{1}{\sqrt{n}} \sum_{i=1}^n \Psi^r(\tilde{X}_i; \mu_n) \\
&= \frac{1}{\sqrt{n}} \sum_{i=1}^n \Psi^r(\tilde{X}_i; \mu) \\
&\quad + \frac{1}{n} \sum_{i=1}^n \sum_{r'=1}^d D_{r'} \Psi^r(\tilde{X}_i; \mu) \sqrt{n}(\mu_n^{r'} - \mu^{r'}) + R_n^r \qquad (1 \le r \le d),
\end{aligned}
$$
(2.7)

where

$$(2.8) \quad R_n^r = \sum_{r'=1}^d \sqrt{n}(\mu_n^{r'} - \mu^{r'}) \frac{1}{n} \sum_{i=1}^n \{D_{r'} \Psi^r(\tilde{X}_i; \theta_n) - D_{r'} \Psi^r(\tilde{X}_i; \mu)\}$$

and $\theta_n$ lies on the line segment joining $\mu$ and $\mu_n$ (for sufficiently large $n$). We will assume

$$(2.9) \quad \begin{aligned} E|\Psi(\tilde{X}_i; \mu)|^2 &< \infty, \\ E|D_{r'} \Psi^r(\tilde{X}_i; \mu)|^2 &< \infty \qquad (\forall\, r, r'). \end{aligned}$$

To show that $R_n^r$ is negligible, write

$$u^{r,r'}(x, \varepsilon) := \sup_{\{\theta \,:\, \|\theta - \mu\| \le \varepsilon\}} |D_{r'} \Psi^r(x; \theta) - D_{r'} \Psi^r(x; \mu)|$$

and assume

$$(2.10) \qquad \delta^{r,r'}(c) := E u^{r,r'}(\tilde{X}_i, c) \to 0 \qquad \text{as } c \downarrow 0 \ (1 \le r, r' \le d).$$

One may then rewrite (2.7) in vectorial form as

$$(2.11) \qquad 0 = \frac{1}{\sqrt{n}} \sum_{i=1}^n \Psi(\tilde{X}_i; \mu) + (\Lambda + \delta_n) \sqrt{n}(\mu_n - \mu),$$

where

$$(2.12) \qquad \Lambda = E((D_{r'} \Psi^r(\tilde{X}_i; \mu)))_{r,r'=1}^d$$

and $\delta_n \to 0$ in probability as $n \to \infty$, if $\mu_n \to \mu$ in probability. If, finally, we assume $\Lambda$ is nonsingular, then (2.11) leads to the equation

$$(2.13) \qquad \sqrt{n}(\mu_n - \mu) = \Lambda^{-1} \left( \frac{1}{\sqrt{n}} \sum_{i=1}^n \Psi(\tilde{X}_i; \mu) \right) + \delta_n',$$

where $\delta_n'$ goes to zero in probability as $n \to \infty$. We have then arrived at the following theorem.

THEOREM 2.1 (CLT for Fréchet sample means). *Let $Q$ be a probability measure on a differentiable manifold $M$ endowed with a metric $\rho$ such that every closed and bounded set of $(M, \rho)$ is compact. Assume* (i) *the Fréchet mean $\mu_{\mathcal{F}}$ exists,* (ii) *there exists a coordinate neighborhood $(U, \phi)$ such that $Q(U) = 1$,* (iii) *the map $\theta \to (\rho^{\phi})^2(\theta, u)$ is twice continuously differentiable on $\phi(U)$,* (iv) *the integrability conditions* (2.9) *hold as well as the relation* (2.10) *and* (v) *$\Lambda$, defined by* (2.12), *is nonsingular. Then* (a) *every measurable selection $\mu_n$ from the (sample) Fréchet mean set of $\hat{Q}_n^{\phi} = \frac{1}{n}\sum_{i=1}^n \delta_{\tilde{X}_i}$ is a consistent estimator of $\mu$, and* (b) $\sqrt{n}(\mu_n - \mu) \xrightarrow{\mathcal{L}} \mathcal{N}(0, \Lambda^{-1} C (\Lambda^t)^{-1})$, *where $C$ is the covariance matrix of $\Psi(\tilde{X}_i; \mu)$.*

PROOF. Part (a) follows from Theorem 2.3 in [8]. The proof of part (b) is as outlined above, and it may also be derived from standard proofs of the CLT for $M$-estimators (see, e.g., [29], pages 132–134). □

As an immediate corollary one obtains:

COROLLARY 2.1. *Let $(M, g)$ be a Riemannian manifold and let $\rho = \rho_g$ be the geodesic distance. Let $Q$ be a probability measure on $M$ whose support is compact and is contained in a coordinate neighborhood $(U, \phi)$. Assume that* (i) *the intrinsic mean $\mu_I = \mu_{\mathcal{F}}$ exists,* (ii) *the map $\theta \to (\rho^{\phi})^2(\theta, u)$ is twice continuously differentiable on $\phi(U)$ for each $u \in \phi(U)$ and $\Lambda$, defined by* (2.12), *is nonsingular. Then the conclusions of Theorem* 2.1 *hold for the intrinsic sample mean $\mu_{n,I} = \mu_{n,\mathcal{F}}$ of $\hat{Q}_n = \frac{1}{n}\sum_{i=1}^n \delta_{X_i}$, with $\mu = \phi(\mu_I)$.*

We now prove one of the main results of this section.

THEOREM 2.2 (CLT for intrinsic sample means). *Let $(M, g)$ be a Riemannian manifold and let $\rho = \rho_g$ be the geodesic distance. Let $Q$ be a probability measure on $M$ whose support is contained in a closed geodesic ball $\overline{B}_r \equiv \overline{B}_r(x_0)$ with center $x_0$ and radius $r$ which is disjoint from the cut locus $C(x_0)$. Assume $r < \frac{\pi}{4K}$, where $K^2$ is the supremum of sectional curvatures in $\overline{B}_r$ if this supremum is positive, or zero if this supremum is nonpositive. Then* (a) *the intrinsic mean $\mu_I$ (of $Q$) exists, and* (b) *the conclusion of Theorem* 2.1 *holds for the image $\mu_n = \phi(\mu_{n,I})$ of the intrinsic sample mean $\mu_{n,I}$ of $\hat{Q}_n = \frac{1}{n}\sum_{i=1}^n \delta_{X_i}$, under the inverse $\phi$ of the exponential map, $\phi = (\mathrm{Exp}_{x_0})^{-1}$.*

PROOF. (a) It is known that under the given assumptions, there is a *local minimum* $\mu_I$, say, of the Fréchet function $F$ which belongs to $B_r$ and that this minimum is also the *unique* minimum in $\overline{B}_{2r}$ [30, 34, 40]. We now show that $\mu_I$ is actually the unique *global minimum* of $F$. Let $p \in (\overline{B}_{2r})^c$. Then $\rho(p,x) > r, \forall x \in \overline{B}_r$. Hence

$$(2.14) \qquad F(p) = \int_{B_r} \rho^2(p,x) Q(dx) > \int_{B_r} r^2 Q(dx) = r^2.$$

On the other hand,

$$(2.15) \qquad F(\mu_I) \leq F(x_0) = \int_{B_r} \rho^2(x_0, x) Q(dx) \leq r^2,$$

proving $F(p) > F(\mu_I)$.

(b) In view of Corollary 2.1, we only need to show that the Hessian matrix $\Lambda \equiv \Lambda(\mu)$ of $F \circ \phi^{-1}$ at $\mu := \phi(\mu_I)$ is nonsingular, where $\phi = \mathrm{Exp}_{x_0}^{-1}$. Now according to [30], Theorem 1.2, for every geodesic curve $\gamma(t)$ in $B_r, t \in (c,d)$ for some $c < 0, d > 0$,

$$(2.16) \qquad \frac{d^2}{dt^2} F(\gamma(t)) > 0 \qquad (c < t < d).$$

Let $\psi = \mathrm{Exp}_{\mu_I}$ denote the exponential map at $\mu_I$, and let $\gamma(t)$ be the unique geodesic with $\gamma(0) = \mu_I$ and $\dot\gamma(0) = v$, so that $\gamma(t) = \psi(tv)$. Here we identify the tangent space $T_{\mu_I} M$ with $\mathbb{R}^d$. Applying (2.16) to this geodesic (at $t = 0$), and writing $G = F \circ \psi$, one has

$$(2.17) \qquad \left.\frac{d^2}{dt^2} F(\psi(tv))\right|_{t=0} = \sum v^i v^j (D_i D_j G)(0) > 0 \qquad (\forall\, v \neq 0),$$

that is, the Hessian of $G$ is positive definite at $0 \in \mathbb{R}^d$. If $x_0 = \mu_I$, this completes the proof of (b).

Next let $x_0 \neq \mu_I$. Now $F \circ \phi^{-1} = G \circ (\psi^{-1} \circ \phi^{-1})$ on a domain that includes $\mu = \phi(\mu_I) \equiv (\mathrm{Exp}_{x_0})^{-1}(\mu_I)$. Write $\psi^{-1} \circ \phi^{-1} = f$. Then in a neighborhood of $\mu$,

$$(2.18) \qquad \begin{aligned} \frac{\partial^2 (G \circ f)}{\partial u^r \, \partial u^{r'}}(u) &= \sum_{j,j'} (D_j D_{j'} G)(f(u)) \frac{\partial f^j}{\partial u^r}(u) \frac{\partial f^{j'}}{\partial u^{r'}}(u) \\ &\quad + \sum_j (D_j G)(f(u)) \frac{\partial^2 f^j}{\partial u^r \, \partial u^{r'}}(u). \end{aligned}$$

The second sum in (2.18) vanishes at $u = \mu$, since $(D_j G)(f(\mu)) = (D_j G)(0) = 0$ as $f(\mu) = \psi^{-1}\phi^{-1}(\mu) = \psi^{-1}(\mu_I) = 0$ is a local minimum of $G$. Also $f$ is a diffeomorphism in a neighborhood of $\mu$. Hence, writing $\Lambda_{r,r'}(\mu)$ as the $(r,r')$ element of $\Lambda(\mu)$,

$$\Lambda_{r,r'}(\mu) = \frac{\partial^2 (F \circ \phi^{-1})}{\partial u^r \, \partial u^{r'}}(\mu) = \sum_{j,j'} (D_j D_{j'} G)(0) \frac{\partial f^j}{\partial u^r}(\mu) \frac{\partial f^{j'}}{\partial u^{r'}}(\mu).$$

This shows, along with (2.17), that $\Lambda = \Lambda(\mu)$ is positive definite. $\square$

REMARK 2.2. If the supremum of the sectional curvatures (of a complete manifold $M$) is nonpositive, and the support of $Q$ is contained in $\overline{B}_r$, then the hypotheses of Theorem 2.2 are satisfied, and the conclusions (a), (b) hold. One may apply this even with $r = \infty$.

REMARK 2.3. The assumptions in Theorem 2.2 on the support of $Q$ for the existence of $\mu_I$ are too restrictive for general applications. But without additional structures they cannot be entirely dispensed with, as is easily shown by letting $Q$ be the uniform distribution on the equator of $S^2$. For the complex projective space $\mathbb{C}P^{d/2}$, $d$ even, necessary and sufficient conditions for the existence of the intrinsic mean $\mu_I$ of an absolutely continuous (w.r.t. the volume measure) $Q$ with radially symmetric density are given in [33, 39].

It may be pointed out that it is the assumption of some symmetry, that is, the invariance of $Q$ under a group of isometries, that often causes the intrinsic mean set to contain more than one element (see, e.g., [8], Proposition 2.2). The next result is, therefore, expected to be more generally applicable than Theorem 2.2.

THEOREM 2.3 (CLT for intrinsic sample means). *Let $Q$ be absolutely continuous w.r.t. the volume measure on a Riemannian manifold $(M, g)$. Assume that* (i) *$\mu_I$ exists,* (ii) *the integrability conditions* (2.9) *hold,* (iii) *the Hessian matrix $\Lambda$ of $F \circ \phi^{-1}$ at $\mu = \phi(\mu_I)$ is nonsingular and* (iv) *the covariance matrix $C$ of $\Psi(\tilde{X}_i; \mu)$ is nonsingular. Then $\sqrt{n}(\mu_n - \mu) \xrightarrow{\mathcal{L}} \mathcal{N}(0, \Gamma)$, where $\Gamma = \Lambda^{-1} C (\Lambda^t)^{-1}$.*

This theorem follows from Theorem 2.1 and Remark 2.1.

In order to obtain a confidence region for $\mu_{\mathcal{F}}$ using the CLT in Theorem 2.1 in the traditional manner, one needs to estimate the covariance matrix $\Gamma = \Lambda^{-1} C (\Lambda^t)^{-1}$. For this one may use proper estimates of $\Lambda$ and $C$, namely,

$$
(2.19) \quad \begin{aligned}
\hat{\Lambda}(\theta) &:= \frac{1}{n} \sum_{i=1}^{n} (\operatorname{Grad} \Psi)(\tilde{X}_i, \mu_n), & \hat{C} &= \operatorname{Cov} \hat{Q}_n^{\phi}, \\
\hat{\Gamma} &:= \hat{\Lambda}^{-1} \hat{C} (\hat{\Lambda}^t)^{-1}, & \hat{\Gamma}^{-1} &= \hat{\Lambda}^t \hat{C}^{-1} \hat{\Lambda}.
\end{aligned}
$$

The following corollary is now immediate. Let $\chi^2_{d,1-\alpha}$ denote the $(1-\alpha)$th quantile of the chi-square distribution with $d$ degrees of freedom.

COROLLARY 2.2. *Under the hypothesis of Theorem 2.1, if $C$ is nonsingular, a confidence region for $\mu_{\mathcal{F}}$ of asymptotic level $1 - \alpha$ is given by $U_{n,\alpha} := \phi^{-1}(D_{n,\alpha})$, where $D_{n,\alpha} = \{v \in \phi(U) : n(\mu_n - v)^t \hat{\Gamma}^{-1} (\mu_n - v) \leq \chi^2_{d,1-\alpha}\}$.*

EXAMPLE 2.2. In the case of the sphere $S^2$ from Example 2.1, it follows that if we consider an arbitrary data point $u = (u^1, u^2)$, and a second point $\theta = \mathrm{Log}_p \lambda = (\theta^1, \theta^2)$, and evaluate the matrix of second-order partial derivatives w.r.t. $\theta^1, \theta^2$ of

$$(2.20) \quad G(u, \theta) = \arccos^2 \left( \cos \|u\| + \frac{\sin \|u\|}{\|u\|}(u^1 \theta^1 + u^2 \theta^2) - \frac{1}{2}\|\theta\|^2 \cos \|u\| \right),$$

then

$$(2.21) \quad \frac{\partial^2 G}{\partial \theta^r \, \partial \theta^s}(u; 0) = \frac{2u^r u^s}{\|u\|^2}\left(1 - \frac{\|u\|}{\tan \|u\|}\right) + \frac{2\delta_{rs}\|u\|}{\tan \|u\|},$$

where $\delta_{rs}$ is the Kronecker symbol and $\|u\|^2 = (u^1)^2 + (u^2)^2$. The matrix $\hat{\Lambda} = (\lambda_{rr'})_{r,r'=1,2}$ has the entries

$$(2.22) \qquad \lambda_{rr'} = \frac{1}{n}\sum_{i=1}^n \frac{\partial^2 G}{\partial \theta^r \, \partial \theta^{r'}}(u_i; 0).$$

Assume $\hat{C}$ is the sample covariance matrix of $u_j, j = 1, \ldots, n$; a large sample confidence region for the intrinsic mean is given by Corollary 2.2 with $\mu_n = 0$.

We now turn to the problem of bootstrapping a confidence region for $\mu_{\mathcal{F}}$. Let $X^*_{i,n}$ be i.i.d. with common distribution $\hat{Q}_n$ (conditionally, given $\{X_i : 1 \leq i \leq n\}$). Write $\tilde{X}^*_{i,n} = \phi(X^*_{i,n}), 1 \leq i \leq n$, and let $\mu^*_n$ be a measurable selection from the Fréchet mean set of $\hat{Q}^{*,\phi}_n := \frac{1}{n}\sum_{i=1}^n \delta_{\tilde{X}^*_{i,n}}$. Let $E^*_{n,\alpha}$ be a subset of $\phi(U)$, such that $P^*(\mu^*_n - \mu_n \in E^*_{n,\alpha}) \to 1 - \alpha$ in probability, where $P^*$ denotes the probability under $\hat{Q}_n$.

COROLLARY 2.3. *In addition to the hypothesis of Theorem* 2.1, *assume $C$ is nonsingular. Then $\phi^{-1}(\{(\mu_n - E^*_{n,\alpha}) \cap \phi(U)\})$ is a confidence region for $\mu_{\mathcal{F}}$ of asymptotic level $(1 - \alpha)$.*

PROOF. One may write (2.7) and (2.8) with $\mu$ and $\mu_n$ replaced by $\mu_n$ and $\mu^*_n$, respectively, also replacing $\tilde{X}_i$ by $\tilde{X}^*_i$ in (2.8). To show that a new version of (2.11) holds with similar replacements (also replacing $\Lambda$ by $\hat{\Lambda}$), with a $\delta^*_n$ (in place of $\delta_n$) going to zero in probability, one may apply Chebyshev's inequality with a first-order absolute moment under $\hat{Q}_n$, proving that $\hat{\Lambda}^* - \hat{\Lambda}$ goes to zero in probability. Here $\hat{\Lambda}^* = \frac{1}{n}\sum_{i=1}^n (\mathrm{Grad}\,\Psi)(\tilde{X}^*_i; \mu^*_n)$. One then arrives at the desired version of (2.7), replacing $\mu_n, \mu, \Lambda, \tilde{X}_i$ by $\mu^*_n, \mu_n, \hat{\Lambda}, \tilde{X}^*_i$, respectively, and with the remainder (corresponding to $\delta'_n$) going to zero in probability. $\square$

REMARK 2.4. In Corollary 2.3, we have considered the so-called *percentile bootstrap* of Hall [25] (also see [17]), which does not require the computation of the standard error $\hat{\Lambda}$. For this as well as for the CLT-based confidence region given by Corollary 2.2, one can show that the coverage error is no more than $O_p(n^{-d/(d+1)})$ or $O(n^{-d/(d+1)})$, as the case may be [4]. One may also use the bootstrap distribution of the *pivotal statistic* $n(\mu_n - \mu)^T \hat{\Gamma}^{-1}(\mu_n - \mu)$ to find $c_{n,\alpha}^*$ such that

$$(2.23) \qquad P^*(n(\mu_n^* - \mu_n)^T \hat{\Gamma}^{*-1}(\mu_n^* - \mu_n) \leq c_{n,\alpha}^*) \simeq 1 - \alpha,$$

to find the confidence region

$$(2.24) \qquad D_{n,\alpha}^* = \{v \in \phi(U) : n(\mu_n - v)^T \hat{\Gamma}^{-1}(\mu_n - v) \leq c_{n,\alpha}^*\}.$$

In particular, if $Q$ has a nonzero absolutely continuous component w.r.t. the volume measure on $M$, then so does $Q^\phi$ w.r.t. the Lebesgue measure on $\phi(U)$ (see [13], page 44). Then assuming (a) $c_{n,\alpha}^*$ is such that the $P^*$-probability in (2.23) equals $1 - \alpha + O_p(n^{-2})$ and (b) some additional smoothness and integrability conditions of the third derivatives of $\Psi$, one can show that the coverage error [i.e., the difference between $1 - \alpha$ and $P(\mu \in D_{n,\alpha}^*)$] is $O_p(n^{-2})$ (see [5, 6, 12, 24, 25]). It follows that the coverage error of the confidence region $\phi^{-1}(D_{n,\alpha}^* \cap \phi(U))$ for $\mu_\mathcal{F}$ is also $O(n^{-2})$. We state one such result precisely.

COROLLARY 2.4 (Bootstrapping the intrinsic sample mean). *Suppose the hypothesis of Theorem 2.3 holds. Then*

$$\sup_{r>0} |P^*(n(\mu_n^* - \mu_n)^T \hat{\Gamma}^{*-1}(\mu_n^* - \mu_n) \leq r)$$
$$- P(n(\mu_n - \mu)^T \hat{\Gamma}^{-1}(\mu_n - \mu) \leq r)| = O_p(n^{-2}),$$

*and the coverage error of the pivotal bootstrap confidence region is* $= O_p(n^{-2})$.

REMARK 2.5. The assumption of absolute continuity of $Q$ in Theorem 2.3 is reasonable for most applications. Indeed this is assumed in most parametric models in directional and shape analysis (see, e.g., [15, 52]).

REMARK 2.6. The results of this section may be extended to the two-sample problem, or to paired samples, in a fairly straightforward manner. For example, in the case of paired observations $(X_i, Y_i), i = 1, \ldots, n$, let $X_i$ have (marginal) distribution $Q$, and intrinsic mean $\mu_I$, and let $Q_2$ and $\nu_I$ be the corresponding quantities for $Y_i$. Let $\phi = \text{Exp}_{x_0}^{-1}$ for some $x_0$, and let $\mu, \nu$ and $\mu_n, \nu_n$ be the images under $\phi$ of the intrinsic population and sample means. Then one arrives at the following [see (2.13)]:

$$(2.25) \qquad \sqrt{n}(\mu_n - \mu) - \sqrt{n}(\nu_n - \nu) \xrightarrow{\mathcal{L}} \mathcal{N}(0, \Gamma),$$

where $\Gamma$ is the covariance matrix of $\Lambda_1^{-1}\Psi(\tilde{X}_i;\mu) - \Lambda_2^{-1}\Psi(\tilde{Y}_i;\nu)$. Here $\Lambda_i$ is the Hessian matrix of $F \circ \phi^{-1}$ for $Q_i$ ($i=1,2$). Assume $\Gamma$ is nonsingular. Then a CLT-based confidence region for $\gamma := \mu - \nu$ is given in terms of $\gamma_n := \mu_n - \nu_n$ by $\{v \in \mathbb{R}^d : n(\gamma_n - v)\hat{\Gamma}^{-1}(\gamma_n - v) \leq \chi^2_{d,1-\alpha}\}$. Alternatively, one may use a bootstrap estimate of the distribution of $\sqrt{n}(\gamma_n - \gamma)$ to derive a confidence region.

In Section 5 we consider two applications of results in this section (and one application of the results in Sections 3 and 4). Application 1 deals with the data from a paleomagnetic study of the possible migration of the Earth's magnetic poles over geological time scales. Here $M = S^2$ and the geodesic distance between two points is the arclength between them measured on the great circle passing through them.

Application 3 analyzes some recent three-dimensional image data on the effect of a (temporary) glaucoma-inducing treatment in 12 Rhesus monkeys. On each animal $k = 4$ carefully chosen landmarks are measured on each eye—the normal eye and the treated eye. For each observation (a set of four points in $\mathbb{R}^3$) the effects of translation, rotation and size are removed to obtain a sample of 12 points on the five-dimensional shape orbifold $\Sigma_3^4$. We use the so-called three-dimensional Bookstein coordinates to label these points (see [15], pages 78–80). In order to apply Theorem 2.3 (i.e., its analog indicated above), a somewhat flat Riemannian structure is chosen so that the necessary assumptions can be verified.

**3. The CLT for extrinsic sample means and confidence regions for the extrinsic mean.** From Theorem 2.1 one may derive a CLT for extrinsic sample means similar to Corollary 2.1. In this section, however, we use another approach which, for extrinsic means, is simpler to apply and generally less restrictive.

Recall that the extrinsic mean $\mu_{j,E}(Q)$ of a nonfocal probability measure $Q$ on a manifold $M$ w.r.t. an embedding $j : M \to \mathbb{R}^k$, when it exists, is given by $\mu_{j,E}(Q) = j^{-1}(P_j(\mu))$, where $\mu$ is the mean of $j(Q)$ and $P_j$ is the projection on $j(M)$ (see [8], Proposition 3.1, e.g.). Often the extrinsic mean will be denoted by $\mu_E(Q)$, or simply $\mu_E$, when $j$ and $Q$ are fixed in a particular context. To ensure the existence of the extrinsic mean set, in this section we will assume that $j(M)$ is closed in $\mathbb{R}^k$.

Assume $(X_1, \ldots, X_n)$ are i.i.d. $M$-valued random objects whose common probability distribution is $Q$, and let $\overline{X}_E := \mu_E(\hat{Q}_n)$ be the *extrinsic sample mean*. Here $\hat{Q}_n = \frac{1}{n}\sum_{j=1}^n \delta_{X_j}$ is the empirical distribution.

A CLT for the extrinsic sample mean on a *submanifold* $M$ of $\mathbb{R}^k$ (with $j$ the inclusion map) was derived by Hendriks and Landsman [27] and, independently, by Patrangenaru [44] by different methods. Differentiable manifolds that are not a priori submanifolds of $\mathbb{R}^k$ arise in new areas of data analysis such as in shape analysis, in high-level image analysis, or in signal and image processing (see, e.g., [15, 16, 22, 31, 32, 33, 42, 51]). These manifolds, known under the names of shape spaces and projective shape spaces, are quotient spaces of submanifolds of $\mathbb{R}^k$ (spaces of orbits of actions of Lie groups), rather than submanifolds of $\mathbb{R}^k$. Our approach is a generalization of the adapted frame method of Patrangenaru [44] to closed embeddings in $\mathbb{R}^k$. This method leads to an appropriate dimension reduction in the CLT and, thereby, reduces computational intensity. This method extends the results of Fisher et al. [19] who considered the case $M = S^d$. We expect that with some effort the results of Hendriks and Landsman [27] may be modified to yield the same result.

Assume $j$ is an embedding of a $d$-dimensional manifold $M$ such that $j(M)$ is closed in $\mathbb{R}^k$, and $Q$ is a $j$-nonfocal probability measure on $M$ such that $j(Q)$ has finite moments of order 2 (or of sufficiently high order as needed). Let $\mu$ and $\Sigma$ be, respectively, the mean and covariance matrix of $j(Q)$ regarded as a probability measure on $\mathbb{R}^k$. Let $\mathcal{F}$ be the set of focal points of $j(M)$, and let $P_j : \mathcal{F}^c \to j(M)$ be the projection on $j(M)$. $P_j$ is differentiable at $\mu$ and has the differentiability class of $j(M)$ around any nonfocal point. In order to evaluate the differential $d_\mu P_j$ we consider a special orthonormal frame field that will ease the computations. Assume $p \to (f_1(x), \ldots, f_d(x))$ is a local frame field on an open subset of $M$ such that, for each $x \in M$, $(dj(f_1(x)), \ldots, dj(f_d(x)))$ are orthonormal vectors in $\mathbb{R}^k$. A local frame field $(e_1(p), e_2(p), \ldots, e_k(p))$ defined on an open neighborhood $U \subseteq \mathbb{R}^k$ is *adapted to the embedding $j$* if it is an orthonormal frame field and $\forall x \in j^{-1}(U)$, $(e_r(j(x)) = d_p j(f_r(x))$, $r = 1, \ldots, d$. Let $e_1, e_2, \ldots, e_k$ be the canonical basis of $\mathbb{R}^k$ and assume $(e_1(p), e_2(p), \ldots, e_k(p))$ is an adapted frame field around $P_j(\mu) = j(\mu_E)$. Then $d_\mu P_j(e_b) \in T_{P_j(\mu)} j(M)$ is a linear combination of $e_1(P_j(\mu)), e_2(P_j(\mu)), \ldots, e_d(P_j(\mu))$:

$$(3.1) \qquad d_\mu P_j(e_b) = \sum (d_\mu P_j(e_b)) \cdot e_a(P_j(\mu)) e_a(P_j(\mu)).$$

By the delta method, $n^{1/2}(P_j(\overline{j(X)}) - P_j(\mu))$ converges weakly to $N(0, \Sigma_\mu)$, where $\overline{j(X)} = \frac{1}{n} \sum_{i=1}^n j(X_i)$ and

$$(3.2) \qquad \begin{aligned} \Sigma_\mu = {} & \left[\sum_{a=1}^d d_\mu P_j(e_b) \cdot e_a(P_j(\mu)) e_a(P_j(\mu))\right]_{b=1,\ldots,k} \Sigma \\ & \times \left[\sum d_\mu P_j(e_b) \cdot e_a(P_j(\mu)) e_a(P_j(\mu))\right]^t_{b=1,\ldots,k}. \end{aligned}$$

Here $\Sigma$ is the covariance matrix of $j(X_1)$ w.r.t. the canonical basis $e_1, e_2, \ldots, e_k$. The asymptotic distribution $N(0, \Sigma_\mu)$ is degenerate and can be regarded as a distribution on $T_{P_j(\mu)} j(M)$, since the range of $d_\mu P_j$ is $T_{P_j(\mu)} j(M)$. Note that

$$d_\mu P_j(e_b) \cdot e_a(P_j(\mu)) = 0 \qquad \text{for } a = d+1, \ldots, k.$$

REMARK 3.1. An asymptotic distribution of the extrinsic sample mean can be obtained as a particular case of Theorem 2.1. The covariance matrix in that theorem depends both on the way the manifold is embedded and on the chart used. We provide below an alternative CLT, which applies to arbitrary embeddings, leads to pivots and is independent of the chart used.

The tangential component $\tan(v)$ of $v \in \mathbb{R}^k$ w.r.t. the basis $e_a(P_j(\mu)) \in T_{P_j(\mu)} j(M), a = 1, \ldots, d$, is given by

$$\tan(v) = (e_1(P_j(\mu))^t v, \ldots, e_d(P_j(\mu))^T v)^t. \tag{3.3}$$

Then the random vector $(d_{\mu_E} j)^{-1}(\tan(P_j(j(\overline{X}))) - P_j(\mu))) = \sum_{a=1}^{d} \overline{X}_j^a f_a$ has the following covariance matrix w.r.t. the basis $f_1(\mu_E), \ldots, f_d(\mu_E)$:

$$\begin{aligned}
\Sigma_{j,E} &= e_a(P_j(\mu))^t \Sigma_\mu e_b(P_j(\mu))_{1 \leq a,b \leq d} \\
&= \left[\sum d_\mu P_j(e_b) \cdot e_a(P_j(\mu))\right]_{a=1,\ldots,d} \Sigma \\
&\quad \times \left[\sum d_\mu P_j(e_b) \cdot e_a(P_j(\mu))\right]^t_{a=1,\ldots,d}.
\end{aligned} \tag{3.4}$$

DEFINITION 3.1. The matrix $\Sigma_{j,E}$ given by (3.4) is the *extrinsic covariance matrix* of the $j$-nonfocal distribution $Q$ (of $X_1$) w.r.t. the basis $f_1(\mu_E), \ldots, f_d(\mu_E)$.

When $j$ is fixed in a specific context, the subscript $j$ in $\Sigma_{j,E}$ will be omitted. If, in addition, rank $\Sigma_\mu = d$, $\Sigma_{j,E}$ is invertible and we define the *j-standardized mean vector*

$$\overline{Z}_{j,n} =: n^{1/2} \Sigma_{j,E}^{-1/2} (\overline{X}_j^1, \ldots, \overline{X}_j^d)^T. \tag{3.5}$$

PROPOSITION 3.1. *Assume $\{X_r\}_{r=1,\ldots,n}$ is a random sample from the $j$-nonfocal distribution $j(Q)$, and let $\mu = E(j(X_1))$ and assume the extrinsic covariance matrix $\Sigma_{j,E}$ of $Q$ is finite. Let $(e_1(p), e_2(p), \ldots, e_k(p))$ be an orthonormal frame field adapted to $j$. Then* (a) *the extrinsic sample mean $\overline{X}_E$ has asymptotically a normal distribution in the tangent space to $M$ at $\mu_E(Q)$ with mean 0 and covariance matrix $n^{-1} \Sigma_{j,E}$, and* (b) *if $\Sigma_{j,E}$ is nonsingular, the $j$-standardized mean vector $\overline{Z}_{j,n}$ given in (3.5) converges weakly to $N(0, I_d)$.*

As a particular case of Proposition 3.1, when $j$ is the inclusion map of a submanifold of $\mathbb{R}^k$, we get the following result for nonfocal distributions on an arbitrary closed submanifold $M$ of $\mathbb{R}^k$:

COROLLARY 3.1. *Assume $M \subseteq \mathbb{R}^k$ is a closed submanifold of $\mathbb{R}^k$. Let $\{X_r\}_{r=1,\ldots,n}$ be a random sample from the nonfocal distribution $Q$ on $M$,* and let $\mu = E(X_1)$ and assume the covariance matrix $\Sigma$ of $j(Q)$ is finite. Let $(e_1(p), e_2(p), \ldots, e_k(p))$ be an orthonormal frame field adapted to $M$. Let $\Sigma_E := \Sigma_{j,E}$, where $j: M \to \mathbb{R}^k$ is the inclusion map. Then (a) $n^{1/2} \tan(j(\overline{X}_E) - j(\mu_E))$ converges weakly to $N(0, \Sigma_E)$, and (b) if $\Sigma$ induces a nonsingular bilinear form on $T_{j(\mu_E)}j(M)$, then $\|\overline{Z}_{j,n}\|^2$ converges weakly to the chi-square distribution $\chi_d^2$.

EXAMPLE 3.1. In the case of a hypersphere in $\mathbb{R}^k$, $j(x) = x$ and $P_j = P_M$. We evaluate the statistic $\|\overline{Z}_{j,n}\|^2 = n\|\Sigma_{j,E}^{-1/2} \tan(P_M(\overline{X}) - P_M(\mu))\|^2$. The projection map is $P_M(x) = x/\|x\|$. $P_M$ has the following property: if $v = cx$, then $d_x P_M(v) = 0$; on the other hand, if the restriction of $d_x P_M$ to the orthocomplement of $\mathbb{R}x$ is a conformal map, that is, if $v \cdot x = 0$, then $d_x P_M(v) = \|x\|^{-1} v$. In particular, if we select the coordinate system such that $x = \|x\| e_k$, then one may take $e_a(P_M(x)) = e_a$, and we get

$$d_x P_M(e_b) \cdot e_a(P_M(x)) = \|x\|^{-1} \delta_{ab} \quad \forall a, b = 1, \ldots, k-1, d_x P_M(e_k) = 0.$$

Since $e_k(P_M(\mu))$ points in the direction of $\mu$, $d_\mu P_M(e_b) \cdot \mu = 0, \forall b = 1, \ldots, k-1$, and we get

(3.6) $\quad \Sigma_E = \|\mu\|^{-2} E([\mathbf{X} \cdot e_a(\mu/\|\mu\|)]_{a=1,\ldots,k-1}[\mathbf{X} \cdot e_a(\mu/\|\mu\|)]_{a=1,\ldots,k-1}^t)$

which is the matrix $G$ in formula (A.1) in [19].

REMARK 3.2. The CLT for extrinsic sample means as stated in Proposition 3.1 or Corollary 3.1 cannot be used to construct confidence regions for extrinsic means, since the population extrinsic covariance matrix is unknown. In order to find a consistent estimator of $\Sigma_{j,E}$, note that $\overline{j(X)}$ is a consistent estimator of $\mu$, $d_{\overline{j(X)}} P_j$ converges in probability to $d_\mu P_j$, and $e_a(P_j(\overline{j(X)}))$ converges in probability to $e_a(P_j(\mu))$ and, further,

$$S_{j,n} = n^{-1} \sum (j(X_r) - \overline{j(X)})(j(X_r) - \overline{j(X)})^t$$

is a consistent estimator of $\Sigma$. It follows that

(3.7) $$\left[\sum_{a=1}^d d_{\overline{j(X)}} P_j(e_b) \cdot e_a(P_j(\overline{j(X)})) e_a(P_j(\overline{j(X)}))\right] S_{j,n}$$
$$\times \left[\sum_{a=1}^d d_{\overline{j(X)}} P_j(e_b) \cdot e_a(P_j(\overline{j(X)})) e_a(P_j(\overline{j(X)}))\right]^t$$

is a consistent estimator of $\Sigma_\mu$, and $\tan_{P_j(\overline{j(X)})} v$ is a consistent estimator of $\tan(v)$.

If we take the components of the bilinear form associated with the matrix (3.7) w.r.t. $e_1(P_j(\overline{j(X)})), e_2(P_j(\overline{j(X)})), \ldots, e_d(P_j(\overline{j(X)}))$, we get a consistent estimator of $\Sigma_{j,E}$ given by

$$
(3.8) \quad \begin{aligned} G(j,X) = & \left[\left[\sum d_{\overline{j(X)}}P_j(e_b) \cdot e_a(P_j(\overline{j(X)}))\right]_{a=1,\ldots,d}\right] \cdot S_{j,n} \\ & \times \left[\left[\sum d_{\overline{j(X)}}P_j(e_b) \cdot e_a(P_j(\overline{j(X)}))\right]_{a=1,\ldots,d}\right]^t, \end{aligned}
$$

and obtain the following results.

THEOREM 3.1. *Assume $j \colon M \to \mathbb{R}^k$ is a closed embedding of $M$ in $\mathbb{R}^k$. Let $\{X_r\}_{r=1,\ldots,n}$ be a random sample from the $j$-nonfocal distribution $Q$, and let $\mu = E(j(X_1))$ and assume $j(X_1)$ has finite second-order moments and the extrinsic covariance matrix $\Sigma_{j,E}$ of $X_1$ is nonsingular. Let $(e_1(p), e_2(p), \ldots, e_k(p))$ be an orthonormal frame field adapted to $j$. If $G(j,X)$ is given by (3.8), then for $n$ large enough $G(j,X)$ is nonsingular (with probability converging to 1) and* (a) *the statistic*

$$(3.9) \qquad n^{1/2} G(j,X)^{-1/2} \tan(P_j(\overline{j(X)}) - P_j(\mu))$$

*converges weakly to $N(0,I_d)$, so that*

$$(3.10) \qquad n\|G(j,X)^{-1/2} \tan(P_j(\overline{j(X)}) - P_j(\mu))\|^2$$

*converges weakly to $\chi_d^2$, and* (b) *the statistic*

$$(3.11) \qquad n^{1/2} G(j,X)^{-1/2} \tan_{P_j(\overline{j(X)})}(P_j(\overline{j(X)}) - P_j(\mu))$$

*converges weakly to $N(0,I_d)$, so that*

$$(3.12) \qquad n\|G(j,X)^{-1/2} \tan_{P_j(\overline{j(X)})}(P_j(\overline{j(X)}) - P_j(\mu))\|^2$$

*converges weakly to $\chi_d^2$.*

COROLLARY 3.2. *Under the hypothesis of Theorem 3.1, a confidence region for $\mu_E$ of asymptotic level $1-\alpha$ is given by* (a) $C_{n,\alpha} := j^{-1}(U_{n,\alpha})$, *where* $U_{n,\alpha} = \{\mu \in j(M) : n\|G(j,X)^{-1/2} \tan(P_j(\overline{j(X)}) - P_j(\mu))\|^2 \le \chi_{d,1-\alpha}^2\}$, *or by* (b) $D_{n,\alpha} := j^{-1}(V_{n,\alpha})$, *where* $V_{n,\alpha} = \{\mu \in j(M) : n\|G(j,X)^{-1/2} \times \tan_{P_j(\overline{j(X)})}(P_j(\overline{j(X)}) - P_j(\mu))\|^2 \le \chi_{d,1-\alpha}^2\}$.

Theorem 3.1 and Corollary 3.2 involve pivotal statistics. The advantages of using pivotal statistics in bootstrapping for confidence regions are well known (see, e.g., [1, 2, 5, 9, 24, 25]).

At this point we recall the steps that one takes to obtain a bootstrapped statistic from a pivotal statistic. If $\{X_r\}_{r=1,\ldots,n}$ is a random sample from the unknown distribution $Q$, and $\{X_r^*\}_{r=1,\ldots,n}$ is a random sample from the empirical $\hat{Q}_n$, conditionally given $\{X_r\}_{r=1,\ldots,n}$, then the statistic

$$T(X,Q) = n\|G(j,X)^{-1/2}\tan(P_j(\overline{j(X)}) - P_j(\mu))\|^2$$

given in Theorem 3.1(a) has the bootstrap analog

$$T(X^*, \hat{Q}_n) = n\|G(j,X^*)^{-1/2}\tan_{P_j(\overline{j(X)})}(P_j(\overline{j(X^*)}) - P_j(\overline{j(X)}))\|^2.$$

Here $G(j, X^*)$ is obtained from $G(j, X)$ by substituting $X_1^*, \ldots, X_n^*$ for $X_1, \ldots, X_n$, and $T(X^*, \hat{Q}_n)$ is obtained from $T(X, Q)$ by substituting $X_1^*, \ldots, X_n^*$ for $X_1, \ldots, X_n$, $\overline{j(X)})$ for $\mu$ and $G(j, X^*)$ for $G(j, X)$.

The same procedure can be used for the vector-valued statistic

$$V(X,Q) = n^{1/2}G(j,X)^{-1/2}\tan(P_j(\overline{j(X)}) - P_j(\mu)),$$

and as a result we get the bootstrapped statistic

$$V^*(X^*, \hat{Q}_n) = n^{1/2}G(j,X^*)^{-1/2}\tan_{P_j(\overline{j(X)})}(P_j(\overline{j(X^*)}) - P_j(\overline{j(X)})).$$

For the rest of this section, we will assume that $j(Q)$, when viewed as a measure on the ambient space $\mathbb{R}^k$, *has finite moments of sufficiently high order*. If $M$ is compact, then this is automatic. In the noncompact case finiteness of moments of order 12, along with an assumption of a nonzero absolutely continuous component, is sufficient to ensure an Edgeworth expansion up to order $O(n^{-2})$ of the pivotal statistic $V(X,Q)$ (see [5, 6, 12, 19, 24]). We then obtain the following results:

THEOREM 3.2. *Let $\{X_r\}_{r=1,\ldots,n}$ be a random sample from the $j$-nonfocal distribution $Q$ which has a nonzero absolutely continuous component w.r.t. the volume measure on $M$ induced by $j$. Let $\mu = E(j(X_1))$ and assume the covariance matrix $\Sigma$ of $j(X_1)$ is defined and the extrinsic covariance matrix $\Sigma_{j,E}$ is nonsingular and let $(e_1(p), e_2(p), \ldots, e_k(p))$ be an orthonormal frame field adapted to $j$. Then the distribution function of*

$$n\|G(j,X)^{-1/2}\tan(P_j(\overline{j(X)}) - P_j(\mu))\|^2$$

*can be approximated by the bootstrap distribution function of*

$$n\|G(j,X^*)^{-1/2}\tan_{P_j(\overline{j(X)})}(P_j(\overline{j(X^*)}) - P_j(\overline{j(X)}))\|^2$$

*with a coverage error $O_p(n^{-2})$.*

One may also use nonpivotal bootstrap confidence regions, especially when $G(j, X)$ is difficult to compute. The result in this case is the following (see [4]).

PROPOSITION 3.2. *Under the hypothesis of Proposition 3.1, the distribution function of $n\|\tan(P_j(\overline{j(X)}) - P_j(\mu))\|^2$ can be approximated uniformly by the bootstrap distribution of*

$$n\|\tan_{P_j(\overline{j(X)})}(P_j(\overline{j(X^*)}) - P_j(\overline{j(X)}))\|^2$$

*to provide a confidence region for $\mu_E$ with a coverage error no more than $O_p(n^{-d/(d+1)})$.*

REMARK 3.3. Note that Corollary 3.2(b) provides a computationally simpler scheme than Corollary 3.2(a) for large sample confidence regions; but for bootstrap confidence regions Theorem 3.2, which is the bootstrap analog of Corollary 3.2(a), yields a simpler method. The corresponding $100(1-\alpha)\%$ confidence region is $C_{n,\alpha}^* := j^{-1}(U_{n,\alpha}^*)$ with $U_{n,\alpha}^*$ given by

$$U_{n,\alpha}^* = \{\mu \in j(M) : n\|G(j,X)^{-1/2}\tan(P_j(\overline{j(X)}) - P_j(\mu))\|^2 \leq c_{1-\alpha}^*\},$$
(3.13)

where $c_{1-\alpha}^*$ is the upper $100(1-\alpha)\%$ point of the values

(3.14) $\quad n\|G(j,X^*)^{-1/2}\tan_{P_j(\overline{j(X)})}(P_j(\overline{j(X^*)}) - P_j(\overline{j(X)}))\|^2$

among the bootstrap resamples. One could also use the bootstrap analog of the confidence region given in Corollary 3.2(b) for which the confidence region is $D_{n,\alpha}^* := j^{-1}(V_{n,\alpha}^*)$ with $V_{n,\alpha}^*$ given by

(3.15) $\quad \begin{aligned}V_{n,\alpha}^* = \{&\mu \in j(M): \\ &n\|G(j,X)^{-1/2}\tan_{P_j(\overline{j(X)})}(P_j(\overline{j(X)}) - P_j(\mu))\|^2 \leq d_{1-\alpha}^*\},\end{aligned}$

where $d_{1-\alpha}^*$ is the upper $100(1-\alpha)\%$ point of the values

(3.16) $\quad n\|G(j,X^*)^{-1/2}\tan_{P_j(\overline{j(X^*)})}(P_j(\overline{j(X^*)}) - P_j(\overline{j(X)}))\|^2$

among the bootstrap resamples. The region given by (3.13)–(3.14) has coverage error $O_p(n^{-2})$.

**4. Asymptotic distributions of sample mean axes, Procrustes mean shapes and extrinsic mean planar projective shapes.** In this section we focus on the asymptotic distribution of sample means in axial data analysis and in planar shape data analysis. The axial space is the $(N-1)$-dimensional real projective space $M = \mathbb{R}P^{N-1}$ which can be identified with the sphere $S^{N-1} = \{x \in \mathbb{R}^N | \|x\|^2 = 1\}$ with antipodal points identified (see, e.g., [41]). If $[x] = \{x, -x\} \in \mathbb{R}P^{N-1}, \|x\| = 1$, the tangent space at $[x]$ can be described as

(4.1) $\quad T_{[x]}\mathbb{R}P^{N-1} = \{([x], v), v \in \mathbb{R}^N | v^t x = 0\}.$

We consider here the general situation when the distribution on $\mathbb{R}P^{N-1}$ may not be concentrated. Note that for $N$ odd, $\mathbb{R}P^{N-1}$ cannot be embedded in $\mathbb{R}^N$, since for any embedding of $\mathbb{R}P^{N-1}$ in $\mathbb{R}^k$ with $N$ odd, the first Stiefel–Whitney class of the normal bundle is not zero ([43], page 51).

The *Veronese–Whitney* embedding is defined for arbitrary $N$ by the formula

$$j([x]) = xx^t, \qquad \|x\| = 1. \tag{4.2}$$

The embedding $j$ maps $\mathbb{R}P^{N-1}$ into a $(\frac{1}{2}N(N+1)-1)$-dimensional Euclidean hypersphere in the space $S(N,\mathbb{R})$ of real $N \times N$ symmetric matrices, where the Euclidean distance $d_0$ between two symmetric matrices is

$$d_0(A, B) = \operatorname{Tr}((A-B)^2).$$

This embedding, which was already used by Watson [52], is preferred over other embeddings in Euclidean spaces because it is *equivariant* (see [35]). This means that the special orthogonal group $\operatorname{SO}(N)$ of orthogonal matrices with determinant $+1$ acts as a group of isometries on $\mathbb{R}P^{N-1}$ with the metric of constant positive curvature; and it also acts on the left on $S_+(N,\mathbb{R})$, the set of nonnegative definite symmetric matrices with real coefficients, by $T \cdot A = TAT^t$. Also, $j(T \cdot [x]) = T \cdot j([x]), \forall T \in \operatorname{SO}(N), \forall [x] \in \mathbb{R}P^{N-1}$.

Note that $j(\mathbb{R}P^{N-1})$ is the set of all nonnegative definite matrices in $S(N,\mathbb{R})$ of rank 1 and trace 1. The following result appears in [8].

PROPOSITION 4.1. (a) *The set $\mathcal{F}$ of the focal points of $j(\mathbb{R}P^{N-1})$ in $S_+(N,\mathbf{R})$ is the set of matrices in $S_+(N,\mathbf{R})$ whose largest eigenvalues are of multiplicity at least 2.* (b) *The projection $P_j : S_+(N,\mathbb{R}) \setminus \mathcal{F} \to j(\mathbb{R}P^{N-1})$ assigns to each nonnegative definite symmetric matrix $A$ with a highest eigenvalue of multiplicity 1, the matrix $j([m])$, where $m(\|m\|=1)$ is an eigenvector of $A$ corresponding to its largest eigenvalue.*

The following result of Prentice [46] is also needed in the sequel.

PROPOSITION 4.2 ([46]). *Assume $[X_r], \|X_r\| = 1, r = 1, \ldots, n$, is a random sample from a $j$-nonfocal probability measure $Q$ on $\mathbb{R}P^{N-1}$. Then the $j$-extrinsic sample covariance matrix $G(j, X)$ is given by*

$$G(j,X)_{ab} = n^{-1}(\eta_N - \eta_a)^{-1}(\eta_N - \eta_b)^{-1} \\ \times \sum_r (m_a \cdot X_r)(m_b \cdot X_r)(m \cdot X_r)^2, \tag{4.3}$$

*where $\eta_a, a = 1, \ldots, N$, are eigenvalues of $K := n^{-1} \sum_{r=1}^n X_r X_r^t$ in increasing order and $m_a, a = 1, \ldots, N$, are corresponding linearly independent unit eigenvectors.*

Here we give a proof of (4.3) based on the equivariance of $j$ to prepare the reader for a similar but more complicated formula of the analogous estimator given later for $\mathbb{C}P^{k-2}$.

Since the map $j$ is equivariant, w.l.o.g. one may assume that $j(\overline{X}_E) = P_j(\overline{j(X)})$ is a diagonal matrix, $\overline{X}_E = [m_N] = [e_N]$ and the other unit eigenvectors of $\overline{j(X)} = D$ are $m_a = e_a, \forall a = 1, \ldots, N-1$. We evaluate $d_D P_j$. Based on this description of $T_{[x]}\mathbb{R}P^{N-1}$, one can select in $T_{P_j(D)}j(\mathbb{R}P^{N-1})$ the orthonormal frame $e_a(P_j(D)) = d_{[e_N]}j(e_a)$. Note that $S(N, \mathbf{R})$ has the orthobasis $F_a^b, b \leq a$, where, for $a < b$, the matrix $F_a^b$ has all entries zero except for those in the positions $(a, b), (b, a)$ that are equal to $2^{-1/2}$; also $F_a^a = j([e_a])$. A straightforward computation shows that if $\eta_a, a = 1, \ldots, N$, are the eigenvalues of $D$ in their increasing order, then $d_D P_j(F_a^b) = 0, \forall b \leq a < N$ and $d_D P_j(F_a^N) = (\eta_N - \eta_a)^{-1} e_a(P_j(D))$; from this equation it follows that, if $\overline{j(X)}$ is a diagonal matrix $D$, then the entry $G(j, X)_{ab}$ is given by

$$G(j, X)_{ab} = n^{-1}(\eta_N - \eta_a)^{-1}(\eta_N - \eta_b)^{-1} \sum_r X_r^a X_r^b (X_r^N)^2.$$

Taking $\overline{j(X)}$ to be a diagonal matrix and $m_a = e_a$, (4.3) follows.

Note that $\mu_{E,j} = [\nu_N]$, where $(\nu_a), a = 1, \ldots, N$, are unit eigenvectors of $E(XX^t) = E(j(Q))$ corresponding to eigenvalues in their increasing order. Let $T([\nu]) = n\|G(j, X)^{-1/2} \tan(P_j(\overline{j(X)}) - P_j(E(j(Q))))\|^2$ be the statistic given by (3.10). We can derive now the following theorem as a special case of Theorem 3.1(a).

THEOREM 4.1. *Assume $j$ is the Veronese–Whitney embedding of $\mathbb{R}P^{N-1}$ and $\{[X_r], \|X_r\| = 1, r = 1, \ldots, n\}$ is a random sample from a $j$-nonfocal probability measure $Q$ on $\mathbb{R}P^{N-1}$ that has a nondegenerate $j$-extrinsic variance. Then $T([\nu])$ is given by*

(4.4) $$T([\nu]) = n\nu^t[(\nu_a)_{a=1,\ldots,N-1}] G(j, X)^{-1}[(\nu_a)_{a=1,\ldots,N-1}]^t \nu,$$

*and, asymptotically, $T([\nu])$ has a $\chi^2_{N-1}$ distribution.*

PROOF. Since $j$ is an isometric embedding and the tangent space $T_{[\nu_N]}\mathbb{R}P^{N-1}$ has the orthobasis $\nu_1, \ldots, \nu_{N-1}$, if we select the first elements of the adapted moving frame in Theorem 3.1 to be $e_a(P_j(\nu_{E,j})) = (d_{[\nu_N]}j)(\nu_a)$, then the $a$th tangential component of $P_j(\overline{j(X)}) - P_j(\nu)$ w.r.t. this basis of $T_{P_j(E(j(Q)))}j(\mathbb{R}P^{N-1})$ equals up to a sign the $a$th component of $m - \nu_N$ w.r.t. the orthobasis $\nu_1, \ldots, \nu_{N-1}$ in $T_{[\nu_N]}\mathbb{R}R^{N-1}$, namely $\nu_a^t m$. The result follows now from Theorem 3.1(a). □

REMARK 4.1. If we apply Theorem 3.1(b) to the embedding $j$, we obtain a similar theorem due to Fisher, Hall, Jing and Wood [19], where $T([\nu])$ is replaced by $T([m])$. Similar asymptotic results can be obtained for the large sample distribution of Procrustes means of planar shapes, as we discuss below. Recall that the planar shape space $M = \sum_2^k$ of an ordered set of $k$ points in $\mathbb{C}$ at least two of which are distinct can be identified in different ways with the complex projective space $\mathbb{C}P^{k-2}$ (see, e.g., [8, 31]). Here we regard $\mathbb{C}P^{k-2}$ as a set of equivalence classes $\mathbb{C}P^{k-2} = S^{2k-3}/S^1$ where $S^{2k-3}$ is the space of complex vectors in $\mathbb{C}^{k-1}$ of norm 1, and the equivalence relation on $S^{2k-3}$ is by multiplication with scalars in $S^1$ (complex numbers of modulus 1). A complex vector $z = (z^1, z^2, \ldots, z^{k-1})$ of norm 1 corresponding to a given configuration of $k$ landmarks, with the identification described in [8], can be displayed in the Euclidean plane (complex line) with the superscripts as labels. If, in addition, $r$ is the largest superscript such that $z^r \neq 0$, then we may assume that $z^r > 0$. Using this representative of the projective point $[z]$ we obtain a unique graphical representation of $[z]$, which will be called the *spherical representation*.

The *Veronese–Whitney* (or simply *Veronese*) *map* is the embedding of $\mathbb{C}P^{k-2}$ in the space of Hermitian matrices $S(k-1, \mathbb{C})$ given in this case by $j([z]) = zz^*$, where, if $z$ is considered as a column vector, $z^*$ is the adjoint of $z$, that is, the conjugate of the transpose of $z$. The Euclidean distance in the space of Hermitian matrices $S(k-1, \mathbb{C})$ is $d_0^2(A, B) = \text{Tr}((A-B)(A-B)^*) = \text{Tr}((A-B)^2)$.

Kendall [31] has shown that the Riemannian metric induced on $j(\mathbb{C}P^{k-2})$ by $d_0$ is a metric of constant holomorphic curvature. The associated Riemannian distance is known as the *Kendall distance* and the full group of isometries on $\mathbb{C}P^{k-2}$ with the Kendall distance is isomorphic to the special unitary group $\text{SU}(k-1)$ of all $(k-1) \times (k-1)$ complex matrices $A$ with $A^*A = I$ and $\det(A) = 1$.

A random variable $X = [Z], \|Z\| = 1$, valued in $\mathbb{C}P^{k-2}$ is $j$-nonfocal if the highest eigenvalue of $E[ZZ^*]$ is simple, and then the extrinsic mean of $X$ is $\mu_{j,E} = [\nu]$, where $\nu \in \mathbb{C}^{k-1}, \|\nu\| = 1$, is an eigenvector corresponding to this eigenvalue (see [8]). The extrinsic sample mean $\overline{[z]}_{j,E}$ of a random sample $[z_r] = [(z_r^1, \ldots, z_r^{k-1})], \|z_r\| = 1, r = 1, \ldots, n$, from such a nonfocal distribution exists with probability converging to 1 as $n \to \infty$, and is the same as that given by

$$\overline{[z]}_{j,E} = [m], \tag{4.5}$$

where $m$ is a highest unit eigenvector of

$$K := n^{-1} \sum_{r=1}^n z_r z_r^*. \tag{4.6}$$

This means that $\overline{[z]}_{j,E}$ is the full Procrustes estimate for parametric families such as Dryden–Mardia distributions or complex Bingham distributions for planar shapes [35, 36]. For this reason, $\mu_{j,E} = [m]$ will be called the *Procrustes* mean of $Q$.

PROPOSITION 4.3. *Assume $X_r = [Z_r], \|Z_r\| = 1, r = 1, \ldots, n$, is a random sample from a $j$-nonfocal probability measure $Q$ with a nondegenerate $j$-extrinsic covariance matrix on $\mathbb{C}P^{k-2}$. Then the $j$-extrinsic sample covariance matrix $G(j, X)$ as a complex matrix has the entries*

$$G(j,X)_{ab} = n^{-1}(\eta_{k-1} - \eta_a)^{-1}(\eta_{k-1} - \eta_b)^{-1} \\ \times \sum_{r=1}^{n}(m_a \cdot Z_r)(m_b \cdot Z_r)^*|m_{k-1} \cdot Z_r|^2. \tag{4.7}$$

The proof is similar to that given for Proposition 4.2 and is based on the equivariance of the Veronese–Whitney map $j$ w.r.t. the actions of $SU(k-1)$ on $\mathbb{C}P^{k-2}$ and on the set $S_+(k-1, \mathbb{C})$ of nonnegative semidefinite self-adjoint $(k-1)$ by $(k-1)$ complex matrices (see [8]). Without loss of generality we may assume that $K$ in (4.6) is given by $K = \text{diag}\{\eta_a\}_{a=1,\ldots,k-1}$ and the largest eigenvalue of $K$ is a simple root of the characteristic polynomial over $\mathbb{C}$, with $m_{k-1} = e_{k-1}$ as a corresponding complex eigenvector of norm 1. The eigenvectors over $\mathbb{R}$ corresponding to the smaller eigenvalues are given by $m_a = e_a, m'_a = ie_a, a = 1, \ldots, k-2$, and yield an orthobasis for $T_{[m_{k-1}]}j(\mathbb{C}P^{k-2})$. For any $z \in S^{2k-1}$ which is orthogonal to $m_{k-1}$ in $\mathbb{C}^{k-1}$ w.r.t. the real scalar product, we define the path $\gamma_z(t) = [\cos t m_{k-1} + \sin t z]$. Then $T_{P_j(K)}j(\mathbb{C}P^{k-2})$ is generated by the vectors tangent to such paths $\gamma_z(t)$ at $t = 0$. Such a vector, as a matrix in $S(k-1, \mathbb{C})$, has the form $zm^*_{k-1} + m_{k-1}z^*$. In particular, since the eigenvectors of $K$ are orthogonal w.r.t. the complex scalar product, one may take $z = m_a, a = 1, \ldots, k-2$, or $z = im_a, a = 1, \ldots, k-2$, and thus get an orthobasis in $T_{P_j(K)}j(M)$. When we norm these vectors to have unit lengths we obtain the orthonormal frame

$$e_a(P_j(K)) = d_{[m_{k-1}]}j(m_a) = 2^{-1/2}(m_a m^*_{k-1} + m_{k-1} m^*_a),$$
$$e'_a(P_j(K)) = d_{[m_{k-1}]}j(im_a) = i2^{-1/2}(m_a m^*_{k-1} - m_{k-1} m^*_a).$$

Since the map $j$ is equivariant we may assume that $K$ is diagonal. In this case $m_a = e_a$, $e_a(P_j(K)) = 2^{-1/2} E_a^{k-1}$ and $e'_a(P_j(K)) = 2^{-1/2} F_a^{k-1}$, where $E_a^b$ has all entries zero except for those in the positions $(a, b)$ and $(b, a)$ that are equal to 1, and $F_a^b$ is a matrix with all entries zero except for those in the positions $(a, b)$ and $(b, a)$ that are equal to $i$, respectively $-i$. Just as in the real case, a straightforward computation shows that $d_K P_j(E_a^b) = d_K P_j(F_a^b) = 0, \forall a \leq b < k-1$, and

$$d_K P_j(E_a^{k-1}) = (\eta_{k-1} - \eta_a)^{-1} e_a(P_j(K)),$$
$$d_K P_j(F_a^{k-1}) = (\eta_{k-1} - \eta_a)^{-1} e'_a(P_j(K)).$$

We evaluate the extrinsic sample covariance matrix $G(j,X)$ given in (3.8) using the real scalar product in $S(k-1,\mathbb{C})$, namely, $U \cdot V = \operatorname{Re}\operatorname{Tr}(UV^*)$. Note that

$$d_K P_j(E_b^{k-1}) \cdot e_a(P_j(K)) = (\eta_{k-1} - \eta_a)^{-1}\delta_{ba},$$
$$d_K P_j(E_b^{k-1}) \cdot e'_a(P_j(K)) = 0$$

and

$$d_K P_j(F_b^{k-1}) \cdot e'_a(P_j(K))^t = (\eta_{k-1} - \eta_a)^{-1}\delta_{ba},$$
$$d_K P_j(F_b^{k-1}) \cdot e_a(P_j(K)) = 0.$$

Thus we may regard $G(j,X)$ as a complex matrix noting that in this case we get

$$(4.8) \quad \begin{aligned} G(j,X)_{ab} &= n^{-1}(\eta_{k-1} - \eta_a)^{-1}(\eta_{k-1} - \eta_b)^{-1} \\ &\quad \times \sum_{r=1}^n (e_a \cdot Z_r)(e_b \cdot Z_r)^* |e_{k-1} \cdot Z_r|^2, \end{aligned}$$

thus proving (4.7) when $K$ is diagonal. The general case follows by equivariance. We consider now the statistic

$$T(\overline{(X)}_E, \mu_E) = n\|G(j,X)^{-1/2}\tan(P_j(\overline{j(X)}) - P_j(\mu_E))\|^2$$

given in Theorem 3.1 in the present context of random variables valued in complex projective spaces to get:

THEOREM 4.2. *Let $X_r = [Z_r]$, $\|Z_r\| = 1$, $r = 1,\ldots,n$, be a random sample from a Veronese-nonfocal probability measure $Q$ on $\mathbb{C}P^{k-2}$. Then the quantity* (3.10) *is given by*

$$(4.9) \quad T([m],[\nu]) = n[(m \cdot \nu_a)_{a=1,\ldots,k-2}]G(j,X)^{-1}[(m \cdot \nu_a)_{a=1,\ldots,k-2}]^*$$

*and asymptotically $T([m],[\nu])$ has a $\chi^2_{2k-4}$ distribution.*

PROOF. The tangent space $T_{[\nu_{k-1}]}\mathbb{C}P^{k-2}$ has the orthobasis $\nu_1,\ldots,\nu_{k-2},\nu_1^*,\ldots,\nu_{k-2}^*$. Note that since $j$ is an isometric embedding, we may select the first elements of the adapted moving frame in Corollary 3.1 to be $e_a(P_j(\mu)) = (d_{[\nu_{k-1}]}j)(\nu_a)$, followed by $e_a^*(P_j(\mu)) = (d_{[\nu_{k-1}]}j)(\nu_a^*)$. Then the $a$th tangential component of $P_j(\overline{j(X)}) - P_j(\mu)$ w.r.t. this basis of $T_{P_j(\mu)}j(\mathbb{C}P^{k-2})$ equals up to a sign the component of $m - \nu_{k-1}$ w.r.t. the orthobasis $\nu_1,\ldots,\nu_{k-2}$ in $T_{[\nu_{k-1}]}\mathbb{C}P^{k-2}$, which is $\nu_a^t m$; and the $a^*$th tangential components are given by $\nu_a^{*t}m$, and together (in complex multiplication) they yield the complex vector $[(m \cdot \nu_a)_{a=1,\ldots,k-2}]$. The claim follows from this and from (4.3), as a particular case of Corollary 3.1. □

We may derive from this the following large sample confidence regions.

COROLLARY 4.1. *Assume $X_r = [Z_r]$, $\|Z_r\| = 1$, $r = 1, \ldots, n$, is a random sample from a $j$-nonfocal probability measure $Q$ on $\mathbb{C}P^{k-2}$. An asymptotic $(1-\alpha)$-confidence region for $\mu_E^j(Q) = [\nu]$ is given by $R_\alpha(\mathbf{X}) = \{[\nu] : T([m], [\nu]) \leq \chi^2_{2k-4,\alpha}\}$, where $T([m], [\nu])$ is given in (4.9). If $Q$ has a nonzero absolutely continuous component w.r.t. the volume measure on $\mathbb{C}P^{k-2}$, then the coverage error of $R_\alpha(\mathbf{X})$ is of order $O(n^{-1})$.*

For small samples the coverage error could be quite large, and a bootstrap analogue of Theorem 4.2 is preferable.

THEOREM 4.3. *Let $j$ be the Veronese embedding of $\mathbb{C}P^{k-2}$, and let $X_r = [Z_r]$, $\|Z_r\| = 1$, $r = 1, \ldots, n$, be a random sample from a $j$-nonfocal distribution $Q$ on $\mathbb{C}P^{k-2}$ having a nonzero absolutely continuous component w.r.t. the volume measure on $\mathbb{C}P^{k-2}$. Assume in addition that the restriction of the covariance matrix of $j(Q)$ to $T_{[\nu]}j(\mathbb{C}P^{k-2})$ is nondegenerate. Let $\mu_E(Q) = [\nu]$ be the extrinsic mean of $Q$. For a resample $\{Z_r^*\}_{r=1,\ldots,n}$ from the sample consider the matrix $K^* := n^{-1} \sum Z_r^* Z_r^{*^*}$. Let $(\eta_a^*)_{a=1,\ldots,k-1}$ be the eigenvalues of $K^*$ in their increasing order, and let $(m_a^*)_{a=1,\ldots,k-1}$ be the corresponding unit complex eigenvectors. Let $G^*(j, X)^*$ be the matrix obtained from $G(j, X)$ by substituting all the entries with $*$-entries. Then the bootstrap distribution function of*

$$T([m]^*, [m]) := n[(m_{k-1}^* \cdot m_a^*)_{a=1,\ldots,k-2}] G^*((j,X)^*)^{-1} [(m_{k-1} \cdot m_a^*)_{a=1,\ldots,k-2}]^*$$

*approximates the true distribution function of $T([m], [\nu])$ given in Theorem 4.2 with an error of order $O_p(n^{-2})$.*

REMARK 4.2. For distributions that are reasonably concentrated one may determine a nonpivotal bootstrap confidence region using Corollary 3.1(a). The chart used here features affine coordinates in $\mathbb{C}P^{k-2}$. Recall that the complex space $\mathbb{C}^{k-2}$ can be embedded in $\mathbb{C}P^{k-2}$, preserving collinearity. Such a standard *affine* embedding, missing only a hyperplane at infinity, is $(z^1, \ldots, z^{k-2}) \to [z^1 : \cdots : z^{k-1} : 1]$. This leads to the notion of *affine* coordinates of a point

$$p = [z^1 : \cdots : z^m : z^{k-1}], \qquad z^{k-1} \neq 0,$$

to be defined as

$$(w^1, w^2, \ldots, w^{k-2}) = \left(\frac{z^1}{z^{k-1}}, \ldots, \frac{z^{k-2}}{z^{k-1}}\right).$$

To simplify the notation the simultaneous confidence intervals used in the next section can be expressed in terms of simultaneous *complex confidence intervals*. If $z = x + iy, w = u + iv, x < u, y < v$, then we define the complex interval $(z, w) = \{c = a + ib | a \in (x, u), b \in (y, v)\}$.

**5. Applications.** In this last section we consider three applications.

APPLICATION 1. Here we consider the data set of $n = 50$ South magnetic pole positions (latitudes and longitudes), determined from a paleomagnetic study of New Caledonian laterities ([20], page 278). As an example of application of Section 2, we give a large sample confidence region for the mean location of the South pole based on this data. The sample points to a nonsymmetric distribution on $S^2$; the extrinsic sample mean and the intrinsic sample mean are given by

$$\overline{X}_E = (0.0105208, 0.199101, 0.979922)^t$$

and, using $\overline{X}_E$ as the initial input of the necessary minimization for constructing $\overline{X}_I$,

$$\overline{X}_I = p = (0.004392, 0.183800, 0.982954)^t.$$

From Examples 2.1 and 2.2, select the orthobasis $e_1(p), e_2(p)$ given in (2.3) and the logarithmic coordinates $u^1, u^2$ w.r.t. this basis in $T_p S^2$ defined in (2.4). Then compute the matrix $\hat{\Lambda}$ given in (2.22), to get, using Corollary 2.2, the following 95% asymptotic confidence region for $\mu_I$:

$$U = \{\operatorname{Exp}_p(u^1 e_1(p) + u^2 e_2(p))|$$
$$16.6786(u^1)^2 - 2.9806 u^1 u^2 + 10.2180(u^1)^2 \le 5.99146\}.$$

Note that Fisher, Lewis and Embleton ([20], page 112) estimate another location parameter, the *spherical median*. The spherical median here refers to the minimizer of the expected geodesic (or, arc) distance to a given point on the sphere. For this paleomagnetism data, their sample median is at $78.9°, 98.4°$, while the extrinsic sample mean is $78.5°, 89.4°$ and the intrinsic sample mean is $79.4°, 88.6°$. These estimates differ substantially from the current position of the South magnetic pole, a difference accounted for by the phenomenon of migration of the Earth's magnetic poles.

APPLICATION 2. As an application of Section 4, we give a nonpivotal bootstrap confidence region for the mean shape of a group of eight landmarks on the skulls of eight-year-old North American children. The sample used is the University School data ([10], pages 400–405). The data set represents coordinates of anatomical landmarks, whose names and position on the skull are given in [10]. The data are displayed in Figure 1. (The presentation of raw data is similar to other known shape data displays such as in [15], page 46.) The shape variable (in our case, shape of the eight landmarks on the upper mid face) is valued in a planar shape space $\mathbb{C}P^6$ (real dimension $= 12$). A spherical representation of a shape in this case consists of seven marked points; in Figure 2 we display a spherical representation of this data set. A representative for the extrinsic sample mean (spherical representation) is

$$(-0.67151 + 0.66823i, 0.76939 + 1.05712i, -1.03159 - 0.15998i,$$
$$-0.57776 - 0.87257i, 0.77871 - 1.36178i,$$
$$-0.17489 + 0.82106i, 1.00000 + 0.00000i).$$

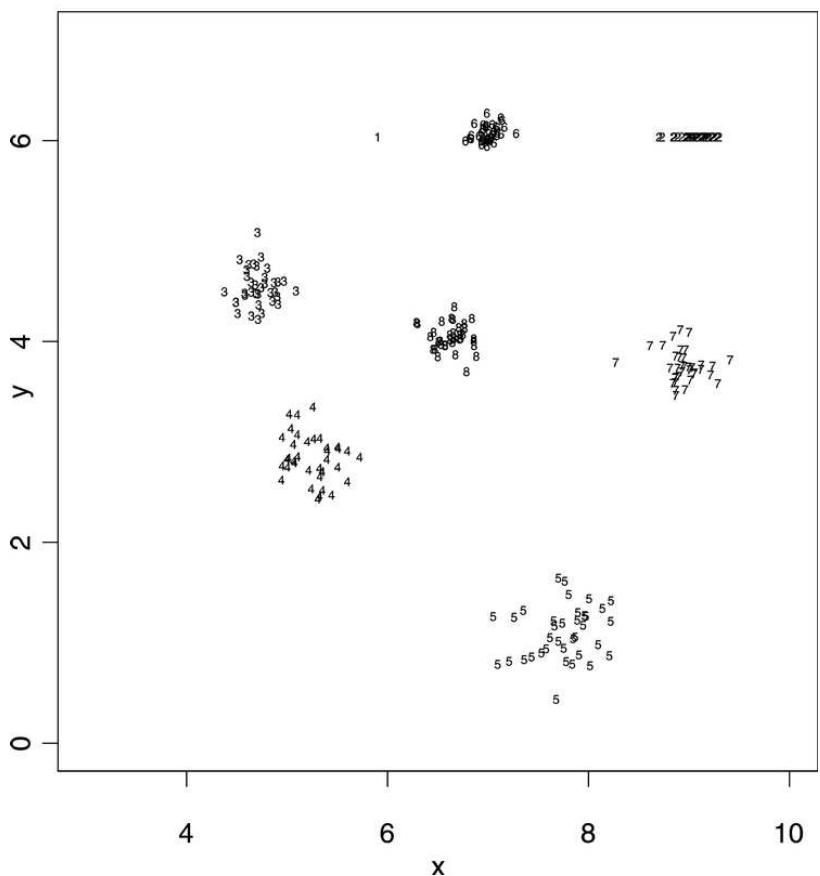

Fig. 1.

We derived the nonpivotal bootstrap distribution using a simple program in S-Plus4.5, that we ran for 500 resamples. A spherical representation of the bootstrap distribution of the extrinsic sample means is displayed in Figure 3. Here we added a representative for the last landmark (the opposite of the sum of the other landmarks since data is centered at 0).

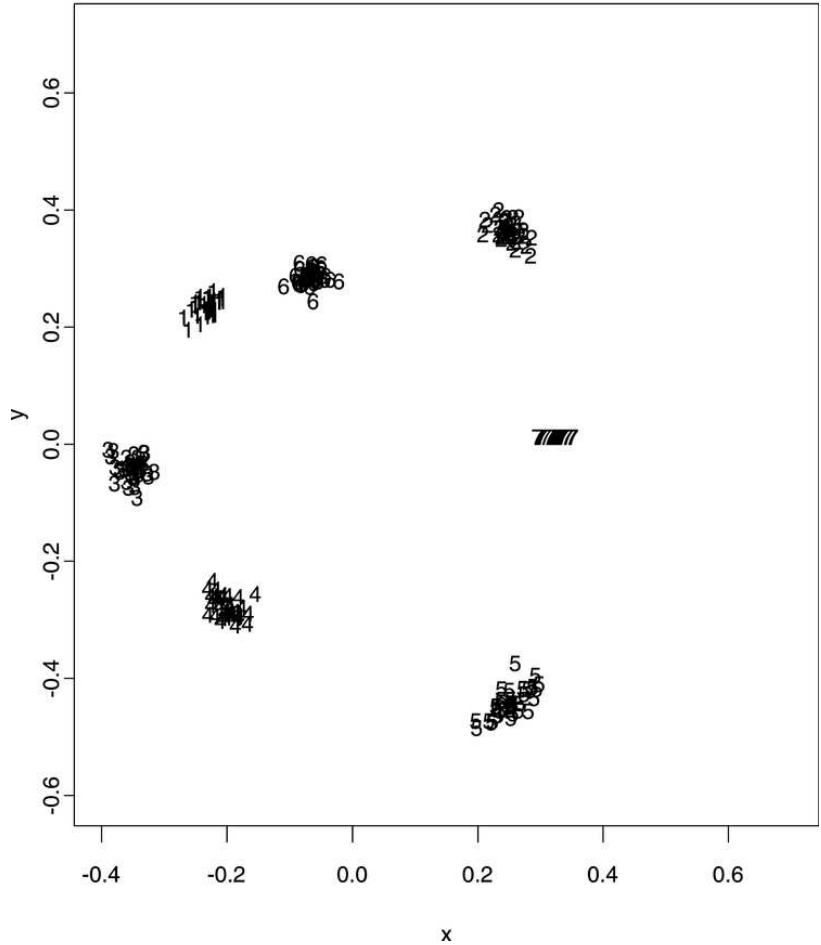

Fig. 2.

Note that the bootstrap distribution of the extrinsic sample mean is very concentrated at each landmark location. This is in agreement with the theory, that predicts in our case a spread of about six times smaller than the spread of the population. It is also an indication of the usefulness of the spherical coordinates. We determined a confidence region for the extrinsic mean using the six 95% simultaneous bootstrap complex intervals for the affine coordinates, as described in Remark 4.2, and found the following complex intervals:

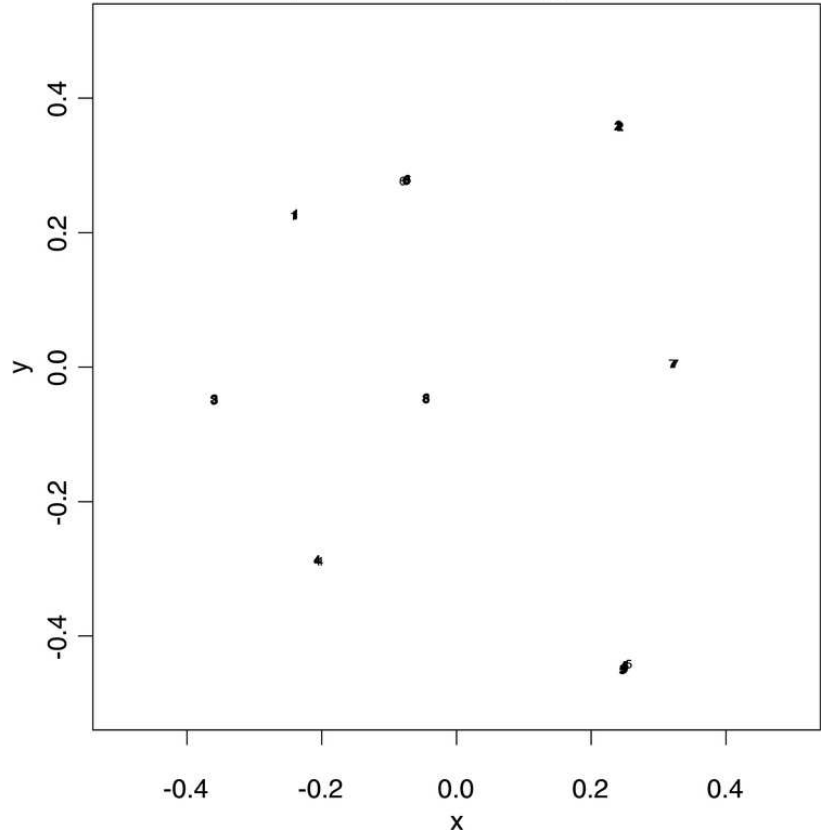

Fig. 3.

for $w_1$:
$$(-0.677268 + 0.666060i, -0.671425 + 0.672409i),$$
for $w_2$:
$$(0.767249 + 1.051660i, 0.775592 + 1.058960i),$$
for $w_3$:
$$(-1.036100 - 0.161467i, -1.029420 - 0.154403i),$$
for $w_4$:
$$(-0.578941 - 0.875168i, -0.574923 - 0.871553i),$$
for $w_5$:
$$(0.777688 - 1.366880i, 0.782354 - 1.358390i),$$

for $w_6$:

$$(-0.177261 + 0.820107i, -0.173465 + 0.824027i).$$

APPLICATION 3. This example is relevant in glaucoma detection. Although it is known that increased intraocular pressure (IOP) may cause a shape change in the eye cup, which is identified with glaucoma, it does not always lead to this shape change. The data analysis presented shows that the device used for measuring the topography of the back of the eye, as reported in [11], is effective in detecting shape change.

We give a nonpivotal bootstrap confidence region for the mean shape change of the eye cup due to IOP. Glaucoma is an eye disorder caused by IOP that is very high. Due to the increased IOP, as the soft spot where the optic nerve enters the eye is pushed backwards, eventually the optic nerve fibers that spread out over the retina to connect to photoreceptors and other retinal neurons can be compressed and damaged. An important diagnostic tool is the ability to detect, in images of the optic nerve head (ONH), increased depth (cupping) of the ONH structures. Two real data-processed images of the ONH cup surface before and after the IOP was increased are shown in Figure 4.

The laser image files are, however, huge-dimensional vectors, and their sizes usually differ. Even if we would restrict the study to a fixed size, there is no direct relationship between the eye cup pictured and the coordinates at a given pixel. A useful data reduction process consists in registration of a number of anatomical landmarks that were identified in each of these images. Assume the position vectors of these landmarks are $X_1, \ldots, X_k, k \geq 4$. Two configurations of landmarks have the same shape if they can be superimposed after a translation, a rotation and a scaling. The shape of the configuration $x = (x_1, \ldots, x_k)$ is labelled $o(x)$ and the space $\Sigma_m^k$ of shapes of configurations of $k$ points in $\mathbb{R}^m$ at least two of which are distinct is the *shape space* introduced by Kendall [31].

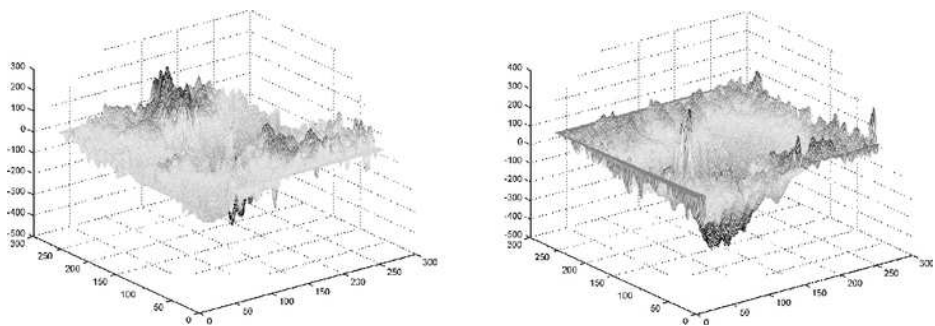

FIG. 4.   *Change in the ONH topography from normal* (left) *to glaucomatous* (right).

We come back to the shape of an ONH cup. This ONH region resembles a "cup" of an ellipsoid and its border has a shape of an ellipse. In this example four landmarks are used. The first three landmarks, denoted by S, T and N, are chosen to be the "top, left and right" points on this ellipse, that is (when referring to the left eye), Superior, Templar and Nose papilla. The last landmark V that we call vertex is the point with the largest "depth" inside the ellipse area that determines the border of the ONH. Therefore, in this example the data analysis is on the shape space of tetrads $\Sigma_3^4$, which is topologically a five-dimensional sphere (see [33], page 38); however, the identification with a sphere is nonstandard. On the other hand, it is known that if a probability distribution on $\Sigma_m^k$ has small support outside a set of singular points, the use of any distance that is compatible with the orbifold topology considered is appropriate in data analysis ([15], page 65) since the data can be linearized. Our choice of the Riemannian metric (5.3) is motivated by considerations of applicability of Theorems 2.2 and 2.3 and computational feasibility. Dryden and Mardia ([15], pages 78–80) have introduced the following five coordinates defined on the generic subset of $\Sigma_3^4$ of shapes of a nondegenerate tetrad that they called *Bookstein coordinates*:

$$
\begin{aligned}
v^1 &= (w_{12}w_{13} + w_{22}w_{23} + w_{32}w_{33})/a, \\
v^2 &= ((w_{12}w_{23} - w_{22}w_{13})^2 \\
    &\qquad + (w_{12}w_{33} - w_{32}w_{13})^2 + (w_{22}w_{33} - w_{23}w_{32})^2)^{1/2}/a, \\
v^3 &= (w_{12}w_{14} + w_{22}w_{24} + w_{32}w_{34})/a, \\
v^4 &= (ab^{1/2})^{-1}(w_{12}^2(w_{23}w_{24} + w_{33}w_{34}) + w_{22}^2(w_{13}w_{14} + w_{33}w_{34}) \\
    &\qquad + w_{32}^2(w_{13}w_{14} + w_{23}w_{24}) - w_{12}w_{13}(w_{22}w_{24} + w_{32}w_{34}) \\
    &\qquad - w_{22}w_{32}(w_{23}w_{34} + w_{33}w_{24}) \\
    &\qquad\qquad\qquad\qquad\qquad - w_{12}w_{14}(w_{22}w_{23} + w_{32}w_{33})), \\
v^5 &= (w_{12}w_{23}w_{34} - w_{12}w_{33}w_{24} - w_{13}w_{22}w_{34} \\
    &\qquad + w_{13}w_{32}w_{24} + w_{22}w_{33}w_{14} - w_{32}w_{23}w_{14})/(2ab)^{1/2},
\end{aligned}
$$
(5.1)

where

$$
\begin{aligned}
a &= 2(w_{12}^2 + w_{22}^2 + w_{32}^2), \\
b &= w_{12}^2 w_{23}^2 + w_{12}^2 w_{33}^2 - 2w_{12}w_{13}w_{22}w_{23} + w_{13}^2 w_{22}^2 + w_{13}^2 w_{32}^2 \\
  &\quad - 2w_{12}w_{13}w_{32}w_{33} + w_{33}^2 w_{22}^2 + w_{23}^2 w_{32}^2 - 2w_{22}w_{32}w_{23}w_{33}
\end{aligned}
$$
(5.2)

and

$$w_{ri} = x_i^r - (x_1^r + x_2^r)/2, \qquad r = 2,3,4.$$

These coordinates carry useful geometric information on the shape of the 4-ad; $v^1$ and $v^3$ give us information on the appearance with respect to the bisector plane of $[X_1 X_2]$, $v^2$ and $v^4$ give some information about the "flatness" of this 4-ad and $v^5$ measures the height of the 4-ad $(X_1, X_2, X_3, X_4)$ relative to the distance $\|X_1 - X_2\|$. Assume $U$ is the set of shapes $o(X)$ such that $(X_1, X_2, X_3, X_4)$ is an affine frame in $\mathbb{R}^3$, and $\phi: U \to \mathbb{R}^{3k-7}$ is the map that associates to $o(X)$ its Bookstein coordinates. $U$ is an open dense set in $\Sigma_3^k$, with the induced topology. In the particular case $k = 4$, $\Sigma_3^4$ is topologically a five-dimensional sphere and, from a classical result of Smale [48], $\Sigma_3^4$ has a differentiable structure diffeomorphic with the sphere $S^5$. Moreover, if $L$ is a compact subset of $U$, there are a finite open covering $U_1 = U, \ldots, U_t$ of $\Sigma_3^4$ and a partition of unity $\phi_1, \ldots, \phi_t$, such that $\phi_1(o(X)) = 1, \forall o(X) \in L$.

We will use the following Riemannian metric on $\Sigma_3^4$: let $(y_1, \ldots, y_5)$ be the Bookstein coordinates of a shape in $U_1$ and let $g_1 = dy_1^2 + \cdots + dy_5^2$ be a flat Riemannian metric on $U_1$, and for each $j = 2, \ldots, t$ we consider any fixed Riemannian metric $g_j$ on $U_j$. Let $g$ be the Riemannian metric given by

$$(5.3) \qquad g = \sum_{j=1}^{t} \phi_j g_j.$$

The space $(\Sigma_3^4, d_g)$ is complete and is flat in a neighborhood of $L$. In this example the two distributions of shapes of tetrads before and after increase in IOP are close. Hence $L$, which contains supports of both distributions, consists of shapes of nondegenerate tetrads only.

Computations for the glaucoma data yield the following results. The $p$-value of the test for equality of the intrinsic means was found to be 0.058, based on the bootstrap distribution of the chi square-like statistic discussed in Remark 2.6. The number of bootstrap resamples for this study was 3000. The chi square-like density histogram is displayed in Figure 5. A matrix plot for the components of the nonpivotal bootstrap distribution of the sample mean differences $\gamma_n^*$ in Remark 2.6 for this application is displayed in Figure 6. The nonpivotal bootstrap 95% confidence intervals for the mean differences $\gamma_j, j = 1, \ldots, 5$, components of $\gamma$ in Remark 2.6 associated with the Bookstein coordinates $v_j, j = 1, \ldots, 5$, are: $(-0.0377073, -0.0058545)$ for $\gamma_1$, $(0.0014153, 0.0119214)$ for $\gamma_2$, $(-0.0303489, 0.0004710)$ for $\gamma_3$, $(0.0031686, 0.0205206)$ for $\gamma_4$, $(-0.0101761, 0.0496181)$ for $\gamma_5$. Note that the individual tests for difference are significant at the 5% level for the first, second and fourth coordinates. However, using the Bonferroni inequality, combining tests for five different shape coordinates each at 5% level leads to a much higher estimated level of significance for the overall shape change.

## APPENDIX

The data set in Application 3 consists of a library of scanning confocal laser tomography (SCLT) images of the complicated ONH topography [11].

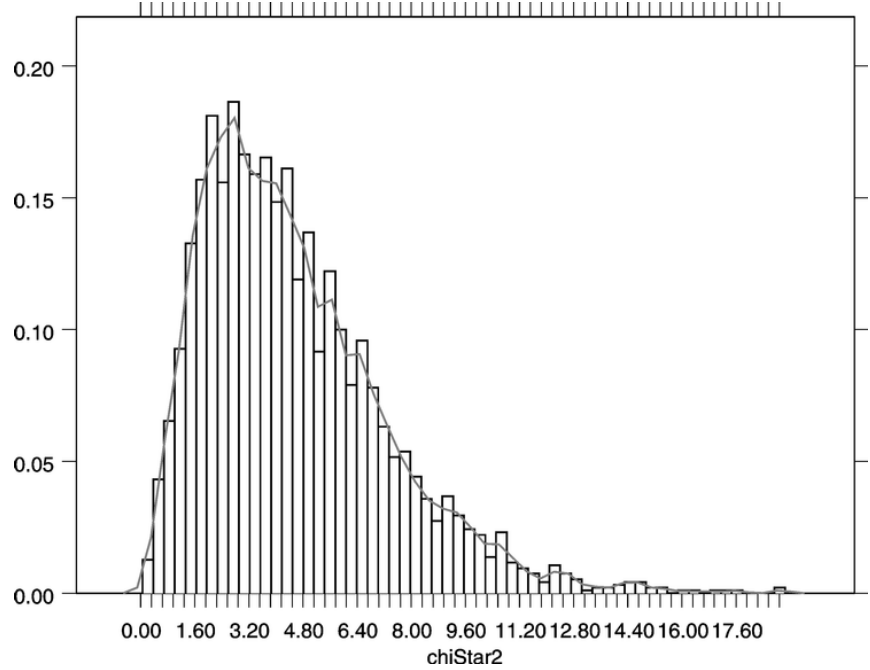

Fig. 5. $\chi^2$-like bootstrap distribution for equality of intrinsic mean shapes from glaucoma data.

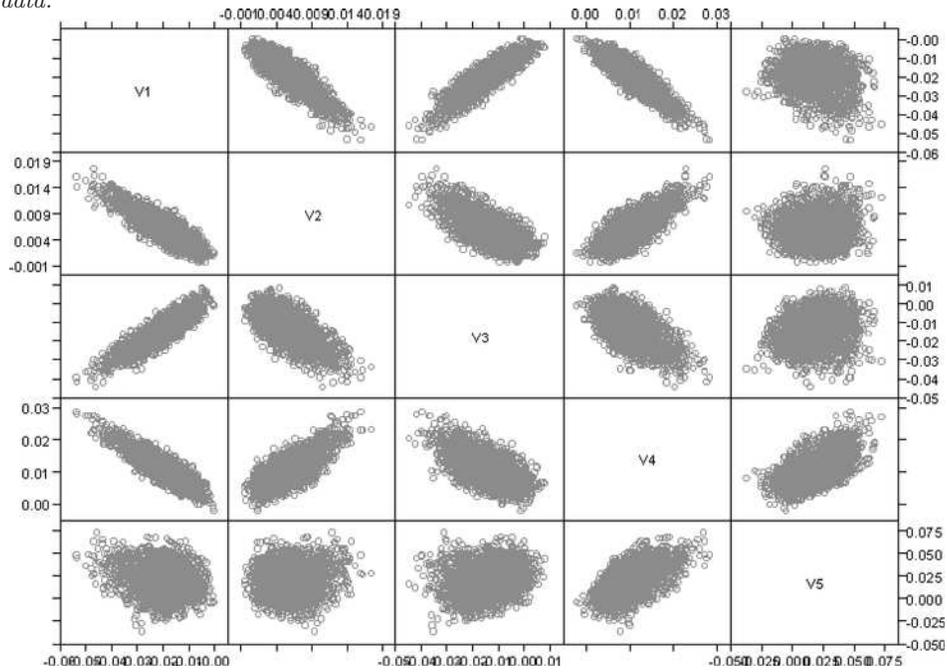

Fig. 6. Glaucoma data, matrix plot for the bootstrap mean differences associated with Bookstein coordinates due to increased IOP.

Those images are the so-called *range images*. A range image is, loosely speaking, like a digital camera image, except that each pixel stores a depth rather than a color level. It can also be seen as a set of points in three dimensions. The range data acquired by 3D digitizers such as optical scanners commonly consist of depths sampled on a regular grid. In the mathematical sense, a range image is a 2D array of real numbers which represent those depths. All of the files (observations) are produced by a combination of modules in C++ and SAS that take the raw image output and process it. The $256 \times 256$ arrays of height values are the products of this software. Another byproduct is a file which we will refer to as the "abxy" file. This file contain the following information: subject names (denoted by: 1c, 1d, 1e, 1f, 1g, 1i, 1j, 1k, 1l, 1n, 1o, 1p), observation points that distinguish the normal and treated eyes and the $10°$ or $15°$ fields of view for the imaging. The observation point "03" denotes a $10°$ view of the experimental glaucoma eye, "04" denotes a $15°$ view of the experimental glaucoma eye, "11" and "12" denote, respectively, the $10°$ and the $15°$ view of the normal eye. The two-dimensional coordinates of the center $(a, b)$ of the ellipses that bound the ONH region, as well as the sizes of the small and the large axes of the ellipses $(x, y)$, are stored in the so-called "abxy" file. To find out more about the LSU study and the image acquisition, see [11]. File names (each file is one observation) were constructed from the information in the "abxy" file. The list of all the observations is then used as an input for the program (created by G. Derado in C++) which determines the three-dimensional coordinates of the landmarks for each observation considered in our analysis, as well as for determining the fifth Bookstein coordinate for each observation. Each image consists of a $256 \times 256$ array of elevation values which represent the "depth" of the ONH. By the "depth" we mean the distance from an imaginary plane, located approximately at the base of the ONH cup, to the "back of the ONH cup."

To reduce the dimensionality of the shape space to 5, out of five landmarks $T$, $S$, $N$, $I$, $V$ recorded, only four landmarks ($X_1 = T$, $X_2 = S$, $X_3 = N$, $X_4 = V$) were considered.

The original data were collected in experimental observations on Rhesus monkeys, and after treatment a healthy eye slowly returns to its original shape. For the purpose of IOP increment detection, in this paper only the first set of after-treatment observations of the treated eye is considered.

**Acknowledgments.** We are thankful to Hilary W. Thompson for providing us with the glaucoma data library and to Sylvie Bogui, Gordana Derado and Jennifer Le for their assistance with part of the computations in Section 5. Thanks are also due to the referees for their thoughtful and constructive suggestions. We also acknowledge the prior work by Hendricks, Landsman and Ruymgaart [28].

DEPARTMENT OF MATHEMATICS
UNIVERSITY OF ARIZONA
TUCSON, ARIZONA 85721
USA
E-MAIL: rabi@math.arizona.edu

DEPARTMENT OF MATHEMATICS
AND STATISTICS
TEXAS TECH UNIVERSITY
LUBBOCK, TEXAS 79402
USA
E-MAIL: vpatrang@math.ttu.edu